\documentclass[11pt, reqno]{amsart}

\usepackage{extarrows}

\usepackage{amsxtra}
\usepackage{hyperref}
\usepackage{tikz-cd}
\usepackage[all]{xy}
\usetikzlibrary{patterns}
\usetikzlibrary{backgrounds}

\let\oldtocsection=\tocsection

\let\oldtocsubsection=\tocsubsection

\let\oldtocsubsubsection=\tocsubsubsection

\renewcommand{\tocsection}[2]{\hspace{0em}\oldtocsection{#1}{#2}}
\renewcommand{\tocsubsection}[2]{\hspace{2em}\oldtocsubsection{#1}{#2}}
\renewcommand{\tocsubsubsection}[2]{\hspace{4.5em}\oldtocsubsubsection{#1}{#2}}

\makeatletter
\def\subsection{\@startsection{subsection}{2}
 \z@{.5\linespacing\@plus.7\linespacing}{-.5em}
 {\normalfont\bfseries}}
\def\subsubsection{\@startsection{subsubsection}{3}
 \z@{.5\linespacing\@plus.7\linespacing}{-.5em}
 {\normalfont\bfseries}}
\makeatother

\theoremstyle{plain}

\newtheorem{theorem}{Theorem}[section]
\newtheorem{lemma}[theorem]{Lemma}
\newtheorem{definition-theorem}[theorem]{Definition-Theorem}
\newtheorem{proposition}[theorem]{Proposition}
\newtheorem{corollary}[theorem]{Corollary}
\newtheorem{definition}[theorem]{Definition}
\newtheorem{example}[theorem]{Example}
\newtheorem{remark}[theorem]{Remark}
\newtheorem{notation}[theorem]{Notation}
\newtheorem{assumption}[theorem]{Assumption}
\newtheorem{lemma-definition}[theorem]{Lemma-Definition}
\newtheorem{lemma-notation}[theorem]{Lemma-Notation}
\newtheorem{question}[theorem]{Question}
\newtheorem{remark-definition}[theorem]{Remarks-Definition}
\newtheorem{notation-remark}[theorem]{Notation-Remarks}
\newtheorem{conjecture}[theorem]{Conjecture}
\newcommand \bth[1] { \begin{theorem}\label{t#1} }
\newcommand \ble[1] { \begin{lemma}\label{l#1} }

\newcommand \bpr[1] { \begin{proposition}\label{p#1} }
\newcommand \bco[1] { \begin{corollary}\label{c#1} }
\newcommand \bconj[1] { \begin{conjecture}\label{co#1} }
\newcommand \bde[1] { \begin{definition}\label{d#1}\rm }
\newcommand \bex[1] { \begin{example}\label{e#1}\rm }
\newcommand \bre[1] { \begin{remark}\label{r#1}\rm }
\newcommand \bnota[1] {\begin{notation}\label{n#1}\rm }
\newcommand \bas[1] { \begin{assumption}\label{a#1}\rm }
\newcommand \bld[1] { \begin{lemma-definition}\label{ld#1} }

\newcommand \bqu[1] { \begin{question}\label{q#1}\rm }

\newcommand {\eth} { \end{theorem} }
\newcommand {\ele} { \end{lemma} }

\newcommand {\epr} { \end{proposition} }
\newcommand {\econj} { \end{conjecture} }
\newcommand {\eco} { \end{corollary} }
\newcommand {\ede} { \end{definition} }
\newcommand {\eex} { \end{example} }
\newcommand {\ere} { \end{remark} }
\newcommand {\enota} { \end{notation} }
\newcommand {\eas} {\end{assumption}}
\newcommand {\eld}{ \end{lemma-definition} }

\newcommand {\equ} {\end{question}}
\newcommand \thref[1]{Theorem \ref{t#1}}
\newcommand \leref[1]{Lemma \ref{l#1}}
\newcommand \prref[1]{Proposition \ref{p#1}}
\newcommand \coref[1]{Corollary \ref{c#1}}
\newcommand \deref[1]{Definition \ref{d#1}}
\newcommand \exref[1]{Example \ref{e#1}}
\newcommand \reref[1]{Remark \ref{r#1}}

\newcommand \notaref[1]{Notation \ref{n#1}}

\newcommand \ldref[1]{Lemma-Definition \ref{ld#1}}

\def \RR {{\mathbb R}}         

\def \ZZ {{\mathbb Z}}

\def \TT {{\mathbb T}}

\def \QQ {{\mathbb Q}}

\def \calQ  {{\mathcal{Q}}}

\def \calP {{\mathcal{P}}}
\def \calX {{\mathcal{X}}}

\renewcommand \max { {\mathrm{max}} }

\newcommand\preceqdot{\mathrel{\ooalign{$\preceq$\cr
  \hidewidth\raise0.225ex\hbox{$\cdot\mkern0.5mu$}\cr}}}

\newcommand{\beqa}{\begin{eqnarray*}}                     
\newcommand{\eeqa}{\end{eqnarray*}}
\def \hs {\hspace{.2in}}

\def \sA {{\scriptscriptstyle A}}

\def \FF {\mathbb{F}}

\def \YY {\mathbb{Y}}
\def \AA {\mathbb{A}}
\def \FF {\mathbb{F}}

\def \calA {{\mathcal{A}}}

\def \calP {{\mathcal{P}}}
\def \calU {{\mathcal{U}}}

\def \mathdash {-\!\!\!-}

\def \AS {{\scriptscriptstyle A,S}}

\def \bfA {{\bf A}}
\def \bfY {{\bf Y}}

\def \calA {\mathcal{A}}
\def \BFZ {{\scriptscriptstyle {\rm BFZ}}}

\def \tkt{{\begin{xy}(0,1)*+{t}="A",(10,1)*+{t'}="B",\ar@{-}^k"A";"B" \end{xy}}}

\def \tit {{\begin{xy}(0,1)*+{t}="A",(10,1)*+{t'}="B",\ar@{-}^i"A";"B" \end{xy}}}
\def \vkv {\begin{xy}(0,1)*+{v}="A",(10,1)*+{v'}="B",\ar@{-}^k"A";"B" \end{xy}}

\begin{document}

\setlength{\baselineskip}{1.2\baselineskip}
\vspace{-.35in}
\title[\bfY-frieze patterns]
{\bfY-frieze patterns}
\author{Antoine de St. Germain}
\address{
Department of Mathematics   \\
The University of Hong Kong \\
Pokfulam Road               \\
Hong Kong}
\email{adsg96@hku.hk}

\date{}

\begin{abstract}
Motivated by cluster ensembles, we introduce a new variant of frieze patterns associated to acyclic cluster algebras, which we call {\bf Y}-frieze patterns. Using the mutation rules for \bfY-variables, we define a large class of \bfY-frieze patterns called {\it unitary \bfY-frieze patterns}, and show that the ensemble map induces a map from (unitary) frieze patterns to (unitary) \bfY-frieze patterns. In rank 2, we show that \bfY-frieze patterns are (associated to) friezes of generalised cluster algebras. In finite type (not necessarily rank 2), we show that \bfY-frieze patterns share the same symmetries as frieze patterns, and prove that their number is finite. 
\end{abstract}

\maketitle
\thispagestyle{empty}
\tableofcontents
\addtocontents{toc}{\protect\setcounter{tocdepth}{1}}

\vspace{1cm}

\section{Introduction and main results}
\subsection{Introduction}\label{ss:intro}
Fix a positive integer $r$, and consider an array  
\begin{equation}\label{eq:frieze-pattern}
r+2\begin{cases} \qquad
    \begin{tikzcd}[column sep=0em, row sep=0.25em]
	\cdots&&1	&&1&&1	&&1&&1&&\cdots&	\\
&& &a_1&&a_2&&a_3&&a_4&&a_5&&\\
	&&\cdots& &b_1&&b_2&&b_3&	&b_4 &&b_5&&\cdots &\\
	&&	&	&&\ddots&&\ddots&&\ddots&&\ddots&&\ddots\\
    &	&&	&\cdots&&1	&&1&&1	&&1&&1 &&\quad\cdots
\end{tikzcd}
\end{cases}
\end{equation}
consisting of $r+2$ staggered infinite rows of positive integers. Recall that such an array is called a {\it Coxeter frieze pattern of width $r$} if every local configuration
\[
\begin{matrix}
			&N&\\
			W&&E\\
			&S&
		\end{matrix}
  \] satisfies the {\it diamond rule}
\[
    WE = 1 + NS.
\]

In this article, we introduce a new variant of Coxeter's frieze patterns. We begin by considering an array  
\begin{equation}\label{eq:Y-frieze-pattern}
r+2\begin{cases} \qquad
    \begin{tikzcd}[column sep=0em, row sep=0.25em]
	\cdots&&0	&&0&&0	&&0&&0&&\cdots&	\\
&& &y_1&&y_2&&y_3&&y_4&&y_5&&\\
	&&\cdots& &z_1&&z_2&&z_3&	&z_4 &&z_5&&\cdots &\\
	&&	&	&&\ddots&&\ddots&&\ddots&&\ddots&&\ddots\\
    &	&&	&\cdots&&0	&&0&&0&&0&&0 &&\quad\cdots
\end{tikzcd}
\end{cases}
\end{equation}
consisting of $r+2$ staggered infinite rows, where the top and bottom rows consist entirely of 0's and the remaining rows consist of positive integers. We say that such an array is a {\it \bfY-frieze pattern of width $r$} if every local configuration
\[
\begin{matrix}
			&N&\\
			W&&E\\
			&S&
		\end{matrix}
  \]
satisfies the {\it \bfY-diamond rule}
\[
    WE = (1+N)(1+S).
\]
     Given a \bfY-frieze pattern as in \eqref{eq:Y-frieze-pattern}, we call the first non-trivial row
\[
    \cdots \quad y_1 \quad y_2 \quad y_3 \quad y_4 \quad y_5 \quad \cdots
\]
 the {\it \bfY-quiddity cycle}. By the \bfY-diamond rule, a \bfY-frieze pattern is determined by its \bfY-quiddity cycle. The following theorem shows that \bfY-frieze patterns and Coxeter frieze patterns are intimately linked. 

 \medskip
\noindent{\bf Theorem A.}
{\it
Each Coxeter frieze pattern of width $r$ as in \eqref{eq:frieze-pattern} gives rise to a unique \bfY-frieze pattern of width $r$ as in \eqref{eq:Y-frieze-pattern} by setting $
y_i = b_i$, for all $i \in \ZZ$.
}
\medskip

 Recall that a {\it glide reflection} is a composition of a reflection and a translation. We denote by ${\bf F}_r$ the unique glide reflection whose square amounts to a horizontal translation by $r+3$. In \cite{Cox}, Coxeter showed that every Coxeter frieze pattern of width $r$ is ${\bf F}_r$-invariant, a property commonly referred to as the {\it glide symmetry} of Coxeter frieze patterns. Remarkably, \bfY-frieze patterns satisfy the same symmetry.
 
\medskip
\noindent{\bf Theorem B.}
{\it
Every \bfY-frieze pattern of width $r$ is ${\bf F}_r$-invariant.}
\medskip

It was recently shown in \cite{muller:finite} that the entries in a frieze pattern of width $r$ admit (explicit) upper bounds, which in turn implies there are finitely many friezes of width $r$. Whilst this is not clear from the definition, similar statements hold for \bfY-frieze patterns. 

\medskip
\noindent{\bf Theorem C.}
{\it
There are finitely many \bfY-frieze patterns of width $r$.}
\medskip

In \cite{Con-Cox1}, Coxeter and Conway constructed a beautiful bijection between Coxeter frieze patterns of width $r$ and triangulations of the regular $(r+3)$-gon. It is a well-known combinatorial fact that triangulations of the $(r+3)$-gon are counted by the Catalan number $C_{r+1}$, where   
\[
    C_r = \frac{1}{r+1} \begin{pmatrix}2r \\ r\end{pmatrix}.
\]

\medskip
\noindent{\bf Conjecture.}
{\it
1) Every \bfY-frieze pattern of width $r$ arises from a Coxeter frieze pattern of width $r$ as in Theorem A. 

2) When $r$ is even, there are $C_{r+1}$-many \bfY-frieze patterns of width $r$. 

3) When $r$ is odd, there are at most $C_{r+1}$-many \bfY-frieze patterns of width $r$.}
\medskip

The elementary nature of Coxeter's frieze patterns conceals a rich algebraic and combinatorial theory related to cluster algebras $\calA$ (in the sense of Fomin and Zelevinsky) of type $A_n$. Unveiling this connection has led to the notion of frieze patterns \cite{ARS:friezes,Keller-Sch:linear-recurrence, Sophie-M:survey}, and their algebraic counterparts called $\calA$-friezes \cite{GM:finite,GHL:friezes}, associated to any acyclic cluster algebra $\calA$. Intuitively, we think of a frieze pattern as the {\it elementary shadow} of an $\calA$-frieze.  

Each (acyclic) cluster algebra $\calA$ naturally belongs to a so-called (acyclic) {\it cluster ensemble} $(\calA, \calX, p)$, where $\calX$ is often called the {\it Poisson cluster algebra} and $p$ the {\it ensemble map}. Our motivation for the notion of \bfY-frieze patterns was a search for (the {\it elementary shadow} of) the ``correct" $\calX$-analog of $\calA$-friezes. The aim of this article is to define \bfY-frieze patterns associated to an arbitrary acyclic cluster ensemble, expose their $\calX$-nature and study their properties in the context of cluster algebras and cluster ensembles.

\subsection{Main results}
In \S \ref{ss:Yfriezepat}, we define frieze patterns and \bfY-frieze patterns associated to an arbitrary symmetrisable generalised Cartan matrix $A$ (\deref{:friezepat}). Since their definitions are subtraction and division-free, we allow frieze patterns and \bfY-frieze patterns to take values in an arbitrary semi-ring $S$ (whose definition is recalled in \S \ref{ss:semi-rings}), thus unifying various notions of patterns appearing in the literature (\reref{:frieze-pattern-history} and \reref{:y-frieze-history}). 

Using the machinery of cluster algebra mutations, S. Morier-Genoud defined in \cite{SMG:arith-2-frieze} a large class of frieze patterns called {\it unitary frieze patterns}. After recalling some notions of mutation in \S \ref{ss:patterns}, we define {\it unitary \bfY-frieze patterns} in \S \ref{ss:unitary}. The following summarises the main results in \S \ref{s:frieze-patterns}.

\medskip
\noindent {\bf Theorem. }(\thref{:frieze-ensemble-map} and \thref{:unitary-to-unitary})
{\it 
Fix an arbitrary symmetrisable generalised Cartan matrix $A$ and a semi-ring $S$. Denote by {\rm Frieze}$(A,S)$, resp. by {\rm YFrieze}$(A,S)$, the set of $S$-valued frieze patterns, resp. \bfY-frieze patterns, associated to $A$. 

1)  One has a well-defined map
\[
    p_\AS : {\rm Frieze}(A,S) \longrightarrow {\rm YFrieze}(A,S).
\]

2) If $f \in {\rm Frieze}(A,S)$ is unitary, then $p_\AS(f) \in {\rm YFrieze}(A,S)$ is {\it unitary}. 
}
\medskip

In \S \ref{s:clusteralg}, we introduce the \bfY-belt algebra (\deref{:belt-alg}), which we use to define the algebraic counterpart of \bfY-frieze patterns (see \deref{:Y-frieze}). It follows rather easily from our setup that \bfY-frieze patterns and their algebraic counterparts are equivalent (\prref{:Y-frieze-pat-to-homo}). This formulation allows us to establish a precise connection between $p_\AS$ and the ensemble map $p$ (\coref{:pas-and-ensemb-map}).

In \S \ref{s:rank2}, we study \bfY-frieze patterns associated to $2 \times 2$ symmetrisable Cartan matrices. The main observation here is that the \bfY-belt algebra associated to a given $A$ is naturally isomorphic to a {\it generalised cluster algebra}, and this isomorphism pulls back (the algebraic counterpart of) \bfY-frieze patterns to friezes of generalised cluster algebras (\prref{:Y-alg-is-gca}). Using this, we determine the number of $\ZZ_{>0}$-valued \bfY-frieze patterns in rank 2 (\thref{:Yfrieze-count-rk2}). 

In \S \ref{s:enum-problem}, we turn to an investigation of \bfY-frieze patterns when the associated Cartan matrix $A$ is of finite type. Under some mild condition on the semi-ring $S$ (see \eqref{eq:embedding-condition}), we determine the symmetry of $S$-valued frieze patterns and $S$-valued \bfY-frieze patterns (\thref{:y-frieze-glide-symmetry}), generalising the well-known {\it glide symmetry} of Coxeter frieze patterns. In \S \ref{ss:enumprob}, we study the problem of enumerating the set of $S$-valued \bfY-frieze patterns when $S$ is the tropical semi-ring (\prref{:friezes-enum}) and the semi-ring of positive integers (\thref{:inf-Yfrieze}).

In this article, \bfY-frieze patterns are associated to an arbitrary Cartan matrix $A$. One deduces Theorems A, B and C in \S \ref{ss:intro} from Theorems \ref{t:frieze-ensemble-map}, \ref{t:y-frieze-glide-symmetry} and \ref{t:inf-Yfrieze} respectively, by taking $A$ of type $A_n$ with the standard labelling (c.f. \S \ref{ss:notation}).

\subsection{Acknowledgements}
 This article began during the author's stay at the MFO in Oberwolfach; we would like to thank the organisers for a wonderful working environment. Preliminary results concerning this work were discussed with S. Morier-Genoud and V. Ovsienko during the author's visit to the University of Champagne-Ardenne. We thank S. Morier-Genoud and V. Ovsienko for their warm welcome, helpful discussions and general enthusiasm for the work presented here. After a first version of this article appeared, we learnt of the proof of \thref{:inf-Yfrieze} from private communication with Greg Muller, whom we thank. We also thank J.-H. Lu for helpful comments on preliminary versions of this article. The author has been partially supported by
the Research Grants Council of the Hong Kong SAR, China (GRF 17307718).

\subsection{Notation}\label{ss:notation}
For any integer $r$, we use $[1, r]$ to denote the set of all integers $k$ from $1$ to $r$.
Elements in $\ZZ^n$ for any integer $n \geq 1$ are regarded as column vectors unless otherwise specified. The standard 
basis vectors in $\ZZ^n$ are denoted by $e_1, e_2, \ldots ,e_n$.  

For a real number $a$, let  $[a]_+={\rm max}(a, 0)$. The transpose of $M$ is denoted by $M^T$.

Suppose that $A$ is an Abelian group. For ${\bf a} = (a_1, \ldots, a_n) \in A^n$ and
$L = (l_1, \ldots, l_n)^T \in \ZZ^n$, we write ${\bf a}^L = a_1^{l_1} \cdots a_n^{l_n} \in A$, and for $L$ an $n \times m$ integral matrix with columns $L_1, \ldots, L_m$, we write ${\bf a}^L = ({\bf a}^{L_1}, \ldots , {\bf a}^{L_m}) \in A^m$. 

Given ${\bf v} = (v_1, \ldots , v_n)$, ${\bf w} = (w_1, \ldots , w_n)$ with $v_i,w_i \in \ZZ$, we write ${\bf v} \leq {\bf w}$ to mean $v_i \leq w_i, i \in [1,n]$.  

{\bf Terminology. } By a {\it symmetrisable generalised Cartan matrix} $A =(a_{i,j})$ we mean a symmetrisable integer matrix such that $a_{i,i} = 2$,  $a_{i,j} \leq 0$ with $a_{i,j} =0$ if and only if $a_{j,i}$ whenever $i \neq j$. A Cartan matrix is of {\it finite type} if all its principal minors are positive, and of {\it infinite type} otherwise. By a Cartan matrix of finite type with the {\it standard labelling}, we mean any Cartan matrix of the form given in \cite[p.43]{Kac:inf-dim-Lie-alg}. 

A Cartan matrix is {\it decomposable} if there exists a non-empty proper subset $I \subset [1,r]$ such that $a_{i,j} = 0$ whenever 
$i \in I$ and $j \in [1,r] \backslash I$. A Cartan matrix which is not decomposable is called {\it indecomposable}.

\section{{\bf Y}-frieze patterns and their properties}\label{s:frieze-patterns}

\subsection{Semi-rings and semi-fields}\label{ss:semi-rings}
 Recall (\cite[\S 2.1]{BFZ:param-canon-basis}) that a {\it semi-ring} $S = (S,\cdot, +)$ (with $1$ and without $0$) consists of a set $S$ and two binary operations $\cdot : S \times S \to S$ (called multiplication) and $+: S \times S \to S$ (called addition) such that 

1) multiplication is associative and commutative, with identity element $1$; 

2) addition is associative and commutative, without $0$;

3) distributivity: $(a+b) c = ac + bc$ for all $a,b,c \in S$.

\bex{:semi-rings}
Our main examples of semi-rings are as follows.

1) The semi-ring $\ZZ_{>0} = (\ZZ_{>0}, \cdot, +)$ of positive integers with the usual addition and multiplication. 

2) The {\it tropical semi-ring} $\ZZ^{\max}_{\geq 0} = (\ZZ_{\geq 0}, +, \max)$ consisting of the set of non-negative integers with semi-ring multiplication the usual addition of integers and semi-ring addition the binary operation $(a,b) \longmapsto \max(a,b)$.
\hfill $\diamond$
\eex

Recall that a semi-ring $(S, \cdot, +)$ is called a {\it semi-field} if $S$ forms a group under multiplication. 

\bex{:semi-fields}
Our main examples of semi-fields are as follows.

    1) The semi-field $\QQ_{>0} = (\QQ_{>0},\cdot,  +)$ of positive rational numbers equipped with the usual addition and multiplication.

    2) The {\it universal semi-field of rank $r$}, denoted $\QQ_{>0}(z_1, \ldots , z_r)$, which consists of the set of non-zero subtraction-free expressions in algebraically independent variables $z_1, \ldots, z_r$ equipped with the operations inherited from the field of rational functions $\QQ(z_1, \ldots, z_r)$. 

    3) The {\it tropical semi-field} $\ZZ^{\max} = (\ZZ, +, \max)$.
    \hfill $\diamond$
\eex
 
\subsection{Frieze patterns and {\bf Y}-frieze patterns}\label{ss:Yfriezepat}
We now introduce the main objects of study of this article. Fix a positive integer $r$.

\bde{:friezepat}
    Let $A = (a_{i,j})$ be any $r \times r$ symmetrisable generalised Cartan matrix, and let $S$ be any semi-ring. 

    1) \cite[Definition 6.2]{GM:finite} An {\it S-valued frieze pattern} associated to $A$ is a map 
    $f :[1,r]\times \ZZ \longrightarrow S$ such that
   \begin{equation}\label{eq:frieze-rel}
       f(i,m) \; f(i,m+1) = 1 + \prod_{j=i+1}^r f(j,m)^{-a_{j,i}} \; \prod_{j=1}^{i-1} f(j,m+1)^{-a_{j,i}}, \qquad \forall \,(i,m) \in [1,r] \times \ZZ.
   \end{equation}

     2) An {\it $S$-valued {\bf Y}-frieze pattern} associated to $A$ is a map $k: [1,r]\times \ZZ \longrightarrow S$ such that
    \begin{equation}\label{eq:Yfrieze-rel}
        k(i,m)\; k(i,m+1) = \prod_{j=i+1}^r (1+k(j,m))^{-a_{j,i}} \; \prod_{j=1}^{i-1}(1+k(j,m+1))^{-a_{j,i}}, \qquad \forall \, (i,m) \in [1,r] \times \ZZ.
    \end{equation}
Let Frieze($A,S$), resp. YFrieze($A,S$), be the set of all $S$-valued frieze patterns, resp. ${\bf Y}$-frieze patterns, associated to $A$. $S$-valued frieze patterns and $S$-valued {\bf Y}-frieze patterns are also said to {\it take values in} $S$. \hfill $\diamond$ 
\ede

\bnota{:infinite-grid}
It is convenient to visualise $[1,r] \times \ZZ$ as an infinite grid of $r$-many staggered rows
\[
    \begin{tikzcd}[column sep=0em, row sep=0em]
    &\cdots&&(1,-1)&&(1,0)&&(1,1)&&(1,2)&&\cdots\\
	&&\cdots&&(2,-1)&&(2,0)&&(2,1)&&(2,2)&&\cdots&\\
	&&&\cdots&&\cdots&&\cdots&&\cdots&&\cdots&& \cdots\\
    &&&&\cdots&&(r,-1)&&(r,0)&&(r,1)&&(r,2)&&\cdots
\end{tikzcd}
\]

An $S$-valued frieze pattern $f$ (resp. an $S$-valued {\bf Y}-frieze pattern $k$) is then an assignment of an element $f(i,m) \in S$ (resp. of an element $k(i,m) \in S$) for each grid point $(i,m)$ of $[1,r]\times \ZZ$, such that \eqref{eq:frieze-rel} (resp. \eqref{eq:Yfrieze-rel}) holds.  \hfill $\diamond$
\enota

\bre{:frieze-pattern-history}
The study of frieze patterns taking values in certain semi-rings has a long history. 
\begin{itemize}
    \item Frieze patterns taking values in $\ZZ_{>0}$ are called {\it arithmetic frieze patterns} in \cite{SMG:arith-2-frieze}, see also \cite{ARS:friezes, GM:finite, Sophie-M:survey}. In particular, when $A$ is of type $A_n$ with the standard labelling, arithmetic frieze patterns are the so-called Coxeter frieze patterns, introduced by Coxeter \cite{Cox} in the 1970s.
    \item Frieze patterns taking values in the tropical semi-field $\ZZ^{\max}$ are called {\it tropical friezes} in \cite{Guo:tropical-frieze}. Tropical friezes were recently studied in the context of the Fock-Goncharov duality conjecture in \cite{CdsgL:ringel-conj}. 
    \item Let $S$ be the universal semi-field of rank $r$, with generators $X_1, \ldots , X_r$. The $S$-valued frieze pattern $F$ associated to $A$ such that
\[
    F(i,0) = X_i, \qquad i \in [1,r],
\]
is called the {\it generic frieze pattern} ({\it generalised frieze pattern} in \cite{Keller-Sch:linear-recurrence}) associated to $A$. 
\end{itemize} 
The notion of $S$-valued frieze patterns given in 1) of \deref{:friezepat} places all the examples above under the same definition. 
\hfill $\diamond$
\ere

\begin{figure}
      \centering
      \[
      \begin{array}{|c|c|}
      \hline
          \text{Cartan type} & \text{Arithmetic {\bf Y}-frieze pattern} \\
          \hline
          A_2 & 
    \begin{tikzcd}[column sep=0em, row sep=0em]
&\cdots& &2&&1&&3&&1&&2&&2&&\cdots\\
	&&\cdots& &1&&2&&2&	&1 &&3&&1&&\cdots
\end{tikzcd}\\
\hline 
G_2 & \begin{tikzcd}[column sep=0em, row sep=0em]
&\cdots& &3&&3&&3&&3&&3&&3&&\cdots\\
	&&\cdots& &8&&8&&8&	&8 &&8&&8&&\cdots
\end{tikzcd} \\
\hline 
A_1^{(1)} & 
\begin{tikzcd}[column sep=0em, row sep=0.25em]
&\cdots&&169&&4&&1&&25&&1156&&\cdots&&&&\\
	&&\cdots&&25&&1&&4&&169&&7921&&\cdots&&&&\\
\end{tikzcd}
\\
\hline
A_3& \begin{tikzcd}[column sep=0em, row sep=0.25em]
&& &2&&2&&2&&2&&2&&2&&\\
	&&\cdots& &3&&3&&3&&3 &&3&&3&&\cdots\\
    &&&& &2&&2&&2&&2&&2&&2&&
\end{tikzcd} \\
\hline
      \end{array}
      \]
      \caption{Arithmetic {\bf Y}-frieze patterns. In type $A_r$, one obtains a {\bf Y}-frieze pattern of width $r$ (c.f. \ref{eq:Y-frieze-pattern}) by adding a top and bottom row of $0$'s.}
      \label{fig:arith-Y-frieze}
  \end{figure}

\bre{:y-frieze-history}
The notion of {\bf Y}-frieze patterns taking values in an arbitrary semi-ring $S$ appears to be new. {\bf Y}-frieze patterns taking values in specific semi-rings have appeared in the literature, but have received far less attention than frieze patterns. 
\begin{itemize}
    \item {\bf Y}-frieze patterns taking values in $\ZZ_{>0}$ are called {\it arithmetic {\bf Y}-frieze patterns}. {\bf Y}-frieze patterns of width $r$ as defined in the Introduction are arithmetic {\bf Y}-frieze patterns associated to a Cartan matrix of type $A_r$ with the standard labelling. Examples of arithmetic {\bf Y}-frieze patterns are given in Figure \ref{fig:arith-Y-frieze}. 
    \item $\ZZ^{\max}$-valued {\bf Y}-frieze patterns are called {\it cluster-additive functions}. They were introduced by Ringel in \cite{Ringel:cluster-additive}. See \cite{CdsgL:ringel-conj} for more.
    \item Let $S$ be the universal semi-field of rank $r$, with generators $Y_1, \ldots , Y_r$. The $S$-valued {\bf Y}-frieze pattern $G$ associated to $A$ such that
\[
    G(i,0) = Y_i, \qquad i \in [1,r]
\]
is called the {\it generic {\bf Y}-frieze pattern} associated to $A$ in \cite[\S 1]{CdsgL:ringel-conj}. Generic {\bf Y}-frieze patterns are a generalisation of the $Y$-systems considered in \cite[\S 8]{FZ:ClusterIV} to study Zamolodchikov's periodicity conjecture. An example of a generic {\bf Y}-frieze pattern is given in Figure \ref{fig:generic-Y-frieze}.
\end{itemize}
\hfill $\diamond$
\ere

Consider the total order on $[1,r] \times \ZZ$ given by
\begin{equation}\label{eq:tot-order-grid}
    (i,m) < (j,n) \quad \iff \quad m < n \quad \text{ or } \quad m= n \text{ and } i < j.
\end{equation}

\ble{:enum-Yfrieze-semifield}
    Let $A$ be any $r \times r$ symmetrisable generalised Cartan matrix, $S$ any semi-field. One has bijections 
    \begin{align*}
        \text{Frieze}(A,S) &\longrightarrow S^r, \qquad f \mapsto (f(1,0), \ldots , f(r,0)) \\
        \text{YFrieze}(A,S) &\longrightarrow S^r, \qquad k \mapsto (k(1,0), \ldots , k(r,0)).
    \end{align*}
    Given ${\bf s} = (s_1, \ldots, s_r) \in S^r$, the unique $f\in {\rm Frieze}(A,S), k \in {\rm YFrieze}(A,S)$ such that 
    \[
        (f(1,0), \ldots , f(r,0)) = {\bf s} \quad \text{and} \quad (k(1,0), \ldots , k(r,0)) = {\bf s},
    \]
    are denoted $f = f_{\bf s}$ and $k = k_{\bf s}$ respectively. The procedure, recursive with respect to the total order \eqref{eq:tot-order-grid}, for computing the values of $f_{\bf s}$ and $k_{\bf s}$ from ${\bf s}$ is called the {\it knitting algorithm} \cite{keller-knitting-algo}.
\ele
\begin{proof}
  Since $S$ is a semi-field, each $f$ is uniquely determined by $(f(1,0), \dots, f(r,0))\in S^r$, by solving \eqref{eq:Yfrieze-rel} recursively using the total order \eqref{eq:tot-order-grid} on $[1,r] \times \ZZ$. It is straightforward to check that two $S$-valued frieze patterns coincide if and only if their restriction to the points $(1,0), \ldots , (r,0)$ coincide. The same argument holds for {\bf Y}-frieze patterns.
\end{proof}

When $S$ is a semi-ring which is not a semi-field, the knitting algorithm fails to produce well-defined $S$-valued maps on $[1,r] \times \ZZ$ in general, since solutions to \eqref{eq:frieze-rel} or \eqref{eq:Yfrieze-rel} in $S$ may not exist. As we shall see in \S \ref{s:enum-problem}, this leads to interesting enumeration problems which depend both on $S$ {\it and} on $A$.

\begin{figure}
    \begin{center}
       \[
    \begin{tikzcd}[column sep=0em, row sep=0.25em]
&& &Y_1&&\frac{Y_2+1}{Y_1}&&\frac{\Delta}{Y_2(1+Y_2)}&&\frac{Y_2+1}{Y_3}\\
	&&\cdots& &Y_2&&\frac{(1+Y_1+Y_2)(1+Y_3)}{Y_1Y_2}&& \frac{(1+Y_1)(1+Y_2+Y_3)}{Y_2Y_3}&&Y_2 &&\cdots\\
    &&&& &Y_3&&\frac{\Delta}{Y_1Y_2Y_3}&&Y_1&&\frac{Y_2 + 1}{Y_1}&&\\
\end{tikzcd}
\] 
    \end{center}
where $\Delta = Y_1Y_2 + Y_1Y_3+Y_2Y_3 + Y_1+Y_2+Y_3+1$. 
   \caption{The generic {\bf Y}-frieze pattern of type $A_3$.}
    \label{fig:generic-Y-frieze}
\end{figure}

\subsection{Constructing {\bf Y}-frieze patterns from frieze patterns}\label{ss:friezepat}
In this subsection, we establish a close connection between frieze patterns and {\bf Y}-frieze patterns. 
\bth{:frieze-ensemble-map}
Let $A$ be any $r\times r$ symmetrisable generalised Cartan matrix and $S$ any semi-ring. One has the set-theoretic map
\[
    p_{\AS} : {\rm Frieze}(A,S) \longrightarrow {\rm YFrieze}(A,S), \qquad f \longmapsto p_{\AS}(f),
\]
given by
\begin{equation}\label{eq:frieze-ensemble-map}
    p_{\AS}(f)(i,m) = \prod_{j=i+1}^r f(j,m)^{-a_{j,i}}\; \prod_{j=1}^{i-1} f(j,m+1)^{-a_{j,i}}.
\end{equation}
\eth
\begin{proof}
    We fix $f \in {\rm Frieze}(A,S)$ and $(i,m) \in [1,r] \times \ZZ$. Set 
    \[
        L = p_{\AS}(f)(i,m) \; p_{\AS}(f)(i,m+1),
    \]
    and
    \[
        R = \prod_{j=i+1}^r (1+p_{\AS}(f)(j,m))^{-a_{j,i}} \; \prod_{j=1}^{i-1} (1+p_{\AS}(f)(j,m+1))^{-a_{j,i}}.
    \]
    We need to show that $p_{\AS}(f)$ satisfies the relation \eqref{eq:Yfrieze-rel}, i.e. that $L = R$.  
    \begin{align*}
        L &= \prod_{j=i+1}^r f(j,m)^{-a_{j,i}} \prod_{j=1}^{i-1} f(j,m+1)^{-a_{j,i}} \; \prod_{j=i+1}^r f(j,m+1)^{-a_{j,i}} \prod_{j=1}^{i-1} f(j,m+2)^{-a_{j,i}}\\
        &= \prod_{j=i+1}^r \left( f(j,m) f(j,m+1)\right)^{-a_{j,i}} \; \prod_{j=1}^{i-1} \left( f(j,m+1) f(j,m+2)\right)^{-a_{j,i}} \\
        &= \prod_{j=i+1}^r \left( 1 + p_{\AS}(f)(j,m)\right)^{-a_{j,i}} \; \prod_{j=1}^{i-1} \left( 1 + p_{\AS}(f)(j,m+1)\right)^{-a_{j,i}} = R,
    \end{align*}
    where in the third equality we used \eqref{eq:frieze-rel}.
\end{proof}

In general, the map $p_{\AS}$ is neither injective nor surjective, as the following example shows.

\bex{:pas-not-bij}
Let $A$ be of type $A_3$ with the standard labelling, let $S = \QQ_{>0}$ and set ${\bf s}= (1,2,3), {\bf t} = (3,2,1), {\bf v} = (1,1,1)$. Consider $ f_{\bf s}, f_{\bf t} \in {\rm Frieze}(A,S)$ and $k_{\bf v} \in {\rm YFrieze}(A,S)$ (see \leref{:enum-Yfrieze-semifield} for the notation). It is readily checked that $p_{\AS}(f_{\bf s}) = p_{\AS}(f_{\bf t})$ and that $k_{\bf v} \not\in {\rm Image}(p_\AS)$. \hfill $\diamond$
\eex

On the other hand, it is easy to find pairs $(A,S)$ for which $p_\AS$ is i) bijective (see \coref{:pas-bij-tropical-sr}); \; ii) injective and not surjective (type $C_2$ with $S = \ZZ_{>0}$); \; iii) surjective and not injective (type $A_3$ with $S = \ZZ_{>0}$). It would be interesting to study this further. 

\subsection{Background on matrix, {\bf A} and \bfY-patterns}\label{ss:patterns} 
In this subsection, we recall various patterns which arise in the theory of cluster algebras, following \cite{CdsgL:ringel-conj}. These will be used in \S \ref{ss:unitary} to define {\it unitary} \bfY-frieze patterns and study their relation to (ordinary) unitary frieze patterns.

A {\it mutation matrix} of size $r$ is any $r \times r$ integer skew-symmetrizable matrix. For a given mutation matrix $B = (b_{i,j})$ of size $r$ and any $k \in [1,r]$, 
the {\it mutation} of $B$ in direction $k$ is the integer matrix $\mu_k(B) = (b_{i,j}')$ where
\begin{equation}\label{eq:matrix-mut}
		b_{i,j}' = \begin{cases}
			-b_{i,j} & \text{ if } i = k \text{ or } j = k, \\
			b_{i,j} + [b_{i,k}]_+[b_{k,j}]_+ - [-b_{i,k}]_+[-b_{k,j}]_+ &\text{ otherwise}.
		\end{cases}
	\end{equation} 

Let $\TT_r$ be an
 $r$-regular tree with the $r$ edges from each vertex labeled bijectively by the set $[1, r]$,
and we write \tkt if the vertices $t$ and $t'$ in $\TT_r$ are
joined by an edge labeled by $k \in [1, r]$. We also fix $t_0 \in \TT_r$ and call it 
 {\it  the root} of $\TT_r$.
 
\bde{:b-pattern}
A {\it matrix pattern of rank $r$} is an assignment, denoted $\{B_t\}_{t \in \TT_r}$, of a mutation matrix $B_t$ of size $r$ to every vertex $t \in \TT_r$ such that $\mu_k(B_t) = B_{t'}$ whenever \tkt.
\hfill $\diamond$

\ede

Let $\FF$ be a field isomorphic to the field of rational functions over $\QQ$ in $r$ independent variables. A {\it seed} in $\FF$ is
a pair $({\bf z}, B)$, where ${\bf z}$ is an ordered set of free generators of $\FF$ over $\QQ$ and $B$ is 
a mutation matrix of size $r$.

Given a seed $({\bf x}, B)$ in $\FF$ and $k \in [1,r]$, the {\it $\bfA$-mutation} of $({\bf x}, B)$ in  direction $k$ is the seed 
$({\bf x}' = (x_1, \dots, x_{k-1}, x_k', x_{k+1}, \dots , x_r), \;\mu_k(B))$,
where $\mu_k(B)$ is given in \eqref{eq:matrix-mut} and
 \begin{equation}\label{eq:X-seed-mutation-cluster}
		x_k' = \frac{\prod_{j=1}^r x_j^{[b_{j,k}]_+} +\prod_{j=1}^r x_j^{[-b_{j,k}]_+}} {x_k}.
\end{equation}

\bde{:cluster-pattern}
1) An {\it $\bfA$-pattern} in $\FF$ is an assignment $\{({\bf x}_t, B_t)\}_{t \in \TT_r}$ 
of a seed $({\bf x}_t, B_t)$ in $\FF$ to each $t \in \TT_r$ such that 
$({\bf x}_{t'}, B_{t'})$ is the $\bfA$-mutation of $({\bf x}_t, B_t)$ in direction $k$ whenever\tkt.

2) The {\it positive space associated to an $\bfA$-pattern $\{({\bf x}_t, B_t)\}_{t \in \TT_r}$} with $B_{t_0} = B$ is the pair
$\AA_B =(\FF(\AA_B), \,\{{\bf x}_t\}_{t \in \TT_r})$. We also call $\AA_B$ {\it the positive $\bfA$-space associated to $B$}.

3) We call each ${\bf x}_t$ a {\it cluster} of $\AA_B$,
each $x \in {\bf x}_t$ a {\it cluster variable} of $\AA_B$, monomials of cluster variables in a single cluster {\it cluster monomials of} $\AA_B$, and $\{B_t\}_{t \in \TT_r}$ the {\it matrix pattern} of $\AA_B$.  
\hfill $\diamond$
\ede

Given a seed $({\bf y}, B)$ in $\FF$ and $k \in [1,r]$, the {\it $\bfY$-mutation} of $({\bf y}, B)$ in 
 direction $k$ is the seed $({\bf y}^\prime, \mu_k(B))$, where 
 $\mu_k(B)$ is given 
 in \eqref{eq:matrix-mut} and  ${\bf y}' = (y_1', \dots, y_r')$ is given by
	\begin{equation}\label{eq:Y-mut}
		y_j' = \begin{cases}
			y_k^{-1}, & \text{ if } j = k; \\
			y_j y_k^{[b_{k,j}]_+} (1+y_k)^{-b_{k,j}}, &\text{ if } j \neq k.
		\end{cases}
	\end{equation}

\bde{:y-pattern}
1) A {\it $\bfY$-pattern} in $\FF$ is an assignment $\{({\bf y}_t, B_t)\}_{t \in \TT_r}$ 
of a seed $({\bf y}_t, B_t)$ in $\FF$ to each $t \in \TT_r$ such that 
$({\bf y}_{t'}, B_{t'})$ is the $\bfY$-mutation of $({\bf y}_t, B_t)$ in direction $k$ whenever\tkt.

2) The {\it positive space associated to a $\bfY$-pattern} $\{({\bf y}_t, B_t)\}_{t \in \TT_r}$ with $B_{t_0} = B$ is the pair $\YY_B = ( \FF(\YY_B), \,\{{\bf y}_t\}_{t \in \TT_r})$. We also call $\YY_B$ {\it the} positive $\bfY$-space associated to $B$.

3) We call each ${\bf y}_t$ a {\it \bfY-cluster} of $\YY_B$,
each $y \in {\bf y}_t$ a {\it $\bfY$-variable of $\YY_B$}, and $\{B_t\}_{t \in \TT_r}$ the {\it matrix pattern} of $\YY_B$. 
\hfill $\diamond$
\ede

Let $B$ be an $r \times r$ mutation matrix, and consider the positive spaces 
\[
\AA_B = (\FF(\AA_B), \; \{{\bf x}_t\}_{t \in \TT_r}) \hs \mbox{and} \hs \YY_B = (\FF(\YY_B), \; \{{\bf y}_t\}_{t \in \TT_r}),
\]
with $\{B_t\}_{t \in \TT_r}$ the matrix pattern of both $\AA_B$ and $\YY_B$. The pair $(\AA_B, \YY_B)$ is called a {\it Fock-Goncharov ensemble}. The sub-algebras
\[
\calU(\AA_B) = \bigcap_{t \in \TT_r} \ZZ[{\bf x}_t^{\pm 1}] \subset \FF(\AA_B) \hs \mbox{and} \hs 
\calU(\YY_B) = \bigcap_{t \in \TT_r} \ZZ[{\bf y}_t^{\pm 1}] \subset \FF(\YY_B)
\]
are called the {\it upper cluster algebra} of $\AA_B$ and of $\YY_B$ respectively. We set
\[
\calU^+(\AA_B) = \bigcap_{t \in \TT_r} \ZZ_{\geq 0}[{\bf x}^{\pm 1}_t] \quad\text{and}\quad \calU^+(\YY_B) = \bigcap_{t \in \TT_r} \ZZ_{\geq 0}[{\bf y}^{\pm 1}_t].
\]

Finally, set $\FF_{>0}(\AA_B) = \QQ_{>0}({\bf x}_t)$ for any $t \in \TT_r$ and $\FF_{>0}(\YY_B) = \QQ_{>0}({\bf y}_t)$ for any $t \in \TT_r$. 
The 
semi-field homomorphism 
\begin{equation}\label{eq:p-map}
    p_{\!{}_B}^*: \FF_{>0}(\YY_B) \longrightarrow  \FF_{>0}(\AA_B), \qquad  
p_{\!{}_B}^*({\bf y}_{t_0}) = {\bf x}_{t_0}^B
\end{equation}

is called the {\it ensemble map} of $(\AA_B, \YY_B)$ in \cite{FG:ensembles}. Note that
$p_{\!{}_B}^*({\bf y}_t) = {\bf x}_t^{B_t}$ for every $t \in \TT_r$.

\subsection{Unitary \bfY-frieze patterns}\label{ss:unitary}

For any $r \times r$ symmetrisable generalised Cartan matrix $A = (a_{i,j})$, consider the acyclic mutation matrix
\begin{equation}\label{eq:B_A}
	B_\sA = \begin{pmatrix}
		0 & a_{1,2} & a_{1,3} & \dots & a_{1,r} \\
		-a_{2,1} & 0 & a_{2,3} & \dots & a_{2,r} \\
		-a_{3,1} & -a_{3,2} & 0 & \dots & a_{3,r} \\
		\dots& \dots & \dots & \dots & \dots \\
		-a_{r,1} & -a_{r,2} & -a_{r,3}& \dots &0
	\end{pmatrix},
\end{equation}
and along with the positive spaces 
\[
    \AA_B = (\FF(\AA_B), \{{\bf x}_t = (x_{t;1}, \ldots , x_{t;r})\}_{t \in \TT_r}) \quad \text{and} \quad \YY_B = (\FF(\YY_B), \{{\bf y}_t = (y_{t;1}, \ldots , y_{t;r})\}_{t \in \TT_r})
\]
associated to $B$. Recall from \cite{GHL:friezes} the 2-regular sub-tree $\TT_r^\flat$
of $\TT_r$ given by
\[ 
\cdots  \stackrel{r-1}{\mathdash} t(r,\!-2) \stackrel{r}{\mathdash} t(1, \!-1) \stackrel{1}{\mathdash} \cdots
\stackrel{r-1}{\mathdash} t(r,\!-1) \stackrel{r}{\mathdash} t(1, \!0) \stackrel{1}{\mathdash} 
\cdots \stackrel{r-1}{\mathdash} t(r,\! 0) \stackrel{r}{\mathdash} t(1,\! 1),
\stackrel{1}{\mathdash} \cdots  
\]
where $t(1, 0) = t_0$. One checks directly that, for every $(i, m) \in [1, r] \times \ZZ$, we have $B_{t(i, m)} = B_{t(i, 0)}$ and that 
\begin{equation}\label{eq:prop-Btim}
    B_{t(i,m)}e_i = (-a_{1,i}, \ldots, -a_{i-1,i}, 0, -a_{i+1,i}, \ldots, -a_{r,i})^T \geq 0.
\end{equation}
where recall that $e_i$ denotes the i$^{\rm th}$ standard basis (column) vector of $\ZZ^r$.

For each $(i,m)$, we write
\[
    x(i,m) \stackrel{\text{def}}{=} x_{t(i,m); i} \in \mathfrak{X}(\AA_B) \quad \text{ and } \quad y(i,m) \stackrel{\text{def}}{=} y_{t(i,m); i} \in \mathfrak{X}(\YY_B).
\]

It is well known that for each $(i,m)$, $x(i,m) \in \calU^+(\AA_B)$.

\ble{:yim-universal-Laurent}\cite[\S 5]{CdsgL:ringel-conj}
For each $(i,m) \in [1,r] \times \ZZ$, we have
\begin{align}
	x(i,m)\, x(i,m+1) &= 1 + \prod_{j=i+1}^r x(j,m)^{-a_{j,i}} \, \,\prod_{j=1}^{i-1} x(j,m+1)^{-a_{j,i}} \label{eq:x(m,i)-relations}\\
 y(i,m)\, y(i,m+1) &= \prod_{j=i+1}^r \left(1 + y(j,m)\right)^{-a_{j,i}} \, \,\prod_{j=1}^{i-1} (1 + y(j,m+1))^{-a_{j,i}},\label{eq:y(m,i)-relations}
\end{align}
and moreover $y(i,m) \in \calU^+(\YY_B)$.
\ele

Given a cluster ${\bf x} = (x_1, \ldots , x_r)$ of $\AA_B$, we denote by $\psi_{\bf x}$ the unique ring homomorphism $\ZZ[{\bf x}^{\pm1}] \to \ZZ$ such that $\psi_{\bf x}(x_i) = 1$ for all $i \in [1,r]$. Similarly, given a \bfY-cluster ${\bf y} = (y_1, \ldots , y_r)$ of $\YY_B$, we denote by $\phi_{\bf y}$ the unique ring homomorphism $\ZZ[{\bf y}^{\pm 1}] \to \ZZ$ such that $\phi_{\bf y}(y_i) = 1$ for all $i \in [1,r]$. 

\bld{:unitary-Y-frieze}
Let $A$ be any symmetrisable generalised Cartan matrix and let $B = B_\sA$ as in \eqref{eq:B_A}. 

1) \cite{SMG:arith-2-frieze} Given a cluster ${\bf x}$ of $\AA_B$, there exists a unique arithmetic frieze pattern $f$ associated to $A$ such that 
\begin{equation}
    f(i,m) = \psi_{\bf x}(x(i,m)), \quad \forall (i,m) \in [1,r] \times \ZZ.
\end{equation}
Arithmetic frieze patterns arising in this way are called {\it unitary}. 

2) Given a \bfY-cluster ${\bf y}$ of $\YY_B$, there exists a unique arithmetic \bfY-frieze pattern $k$ associated to $A$ such that
\begin{equation}
    k(i,m) = \phi_{\bf y}(y(i,m)), \quad \forall (i,m) \in [1,r] \times \ZZ.
\end{equation}
Arithmetic \bfY-frieze patterns obtained in this way are called {\it unitary}.
\eld
\begin{proof}
    This follows directly from \leref{:yim-universal-Laurent}. 
\end{proof}

A subset $\sigma = \{\sigma_1, \ldots , \sigma_r\}$ of $\mathfrak{X}(\AA_B)$ (resp. of $\mathfrak{X}(\YY_B)$) is called an {\it unordered cluster of $\AA_B$} (resp. an {\it unordered \bfY-cluster of $\YY_B$}) if i) $\sigma$ is a generating set of $\FF(\AA_B)$ (resp. of $\FF(\YY_B)$) over $\QQ$; \; ii) there exists a cluster ${\bf x}$ of $\AA_B$ (resp. a \bfY-cluster ${\bf y}$ of $\YY_B$) such that $\sigma_i \in {\bf x}$ (resp. $\sigma_i \in {\bf y}$) for all $i \in [1,r]$. By \ldref{:unitary-Y-frieze}, we have well-defined assignments 
\begin{align}
    \{\text{unordered clusters of } \AA_B\} &\longrightarrow {\rm Frieze}(A,\ZZ_{>0}) \label{eq:uc-to-frieze} \\
     \{\text{unordered \bfY-clusters of } \YY_B\} &\longrightarrow {\rm YFrieze}(A,\ZZ_{>0}) \label{eq:uc-to-Yfrieze} 
\end{align}
It is by now well-known (see e.g. \cite[Proposition 3.17]{GM:finite}) that \eqref{eq:uc-to-frieze} is injective, but in general not surjective (e.g. in type $G_2$). 

On the other hand, the assignment \eqref{eq:uc-to-Yfrieze} is not injective in general. This can be seen already in type $A_1$: there are exactly two \bfY-clusters, both giving rise to the same arithmetic \bfY-frieze pattern. Furthermore, \eqref{eq:uc-to-Yfrieze} is in general not surjective either. This is the case in type $G_2$, where there are 8 unordered \bfY-clusters but 21 arithmetic \bfY-frieze patterns (see \thref{:Yfrieze-count-rk2}).

\bth{:unitary-to-unitary}
Let $A$ be any symmetrisable generalised Cartan matrix and let $S = \ZZ_{>0}$. If $f \in {\rm Frieze}(A, S)$ is unitary, then $p_\AS(f) \in {\rm YFrieze}(A, S)$ is unitary. 
\eth
\begin{proof}
    Let $f \in {\rm Frieze}(A,\ZZ_{>0})$, and suppose that $f$ is unitary. By definition, there exists a cluster ${\bf x}_t = (x_{t;1}, \ldots , x_{t;r}), t \in \TT_r$ of $\AA_B$ such that $f(i,m) = \psi_{{\bf x}_t}(x(i,m))$ for each $(i,m) \in [1,r] \times \ZZ$. Following \cite[Lemma 5.4]{GHL:friezes}, one can show that 
    \begin{equation}\label{eq:xt(i,m)}
        x_{t(i,m)} = (x(1, m+1), \ldots , x(i-1,m+1),x(i,m), x(i+1,m), \ldots , x(r,m)). 
    \end{equation}
    Combining \eqref{eq:xt(i,m)} with \eqref{eq:prop-Btim}, we deduce that
    \begin{equation}\label{eq:LHS-unitary}
        p_{\AS}(f)(i,m) = \psi_{{\bf x}_t}({\bf x}_{t(i,m)})^{B_{t(i,m)}e_i}.
    \end{equation}
    Denote by $\Psi_t: \FF_{>0}(\AA_B) \longrightarrow \QQ_{>0}$ and $\Phi_t: \FF_{>0}(\YY_B) \longrightarrow \QQ_{>0}$ the semi-field homomorphisms determined by $\Psi_t(x_{t;i}) = 1$ and $\Phi_t(y_{t;i}) = 1$ respectively, for $i \in [1,r]$. Set $R_t = \FF_{>0}(\AA_B) \cap \ZZ[{\bf x}_t^{\pm 1}]$ and $R_t' = \FF_{>0}(\YY_B) \cap \ZZ[{\bf y}_t^{\pm 1}]$. Note that $R_t, R_t'$ are closed under multiplication, $x(i,m) \in R_t$ by the Laurent phenomenon \cite{FZ:ClusterIV}, $y(i,m) \in R_t'$ by \leref{:yim-universal-Laurent} and 
    \begin{equation}\label{eq:semifield-ring-restrict}
        \Psi_t\big\vert_{R_t} = \psi_{{\bf x}_t}\big\vert_{R_t}, \qquad \Phi_t\big\vert_{R_t'} = \phi_{{\bf y}_t}\big\vert_{R_t'} \qquad \text{ and } \qquad \Psi_t\circ p^*_{\!{}_B} = \Phi_t.
    \end{equation}
    By applying \eqref{eq:semifield-ring-restrict} to \eqref{eq:LHS-unitary}, we obtain
    \[
        p_{\AS}(f)(i,m) = \Psi_t({\bf x}_{t(i,m)}^{B_{t(i,m)}e_i}) = \Psi_t \circ p^*_{\!{}_B}(y(i,m)) = \Phi_t(y(i,m)) = \phi_{{\bf y}_t}(y(i,m)).
    \]
\end{proof}

We conclude this section with an extended example in type $A_3$.

\bex{:A3-unitary-Y-frieze}
Let $A$ be of type $A_3$ with the standard labelling, set $B = B_\sA$ as in \eqref{eq:B_A}, and consider the positive space 
\[
    \YY_B = (\FF(\YY_B), \{{\bf y}_t = (y_{t;1}, y_{t;2}, y_{t;3})\}_{t \in \TT_3}).
\]
At $t_0 \in \TT_3$, we write ${\bf y}_{t_0} = (y_1,y_2,y_3)$. Note that $y(i,m) \in \ZZ_{\geq 0}[y_1^{\pm 1}, y_2^{\pm 1}, y_3^{\pm 1}]$ by \leref{:yim-universal-Laurent}; see Figure \ref{fig:yim-laurent} for the Laurent expansion of a few $y(i,m)$'s (in fact, one can show that Figure \ref{fig:yim-laurent} contains all the $y(i,m)$'s; see \S\ref{ss:glide-sym}).

\begingroup
\begin{figure}
    \centering
\renewcommand{\arraystretch}{2}
\[
\begin{array}{|c||c|c|c|c|c|}
    \hline
    m & 0 & 1 & 2 & 3 & 4\\
     \hline \hline
    y(1,m) & y_1 & \frac{1 + (1 + y_1)y_2}{y_1} & \frac{1 + (1 + y_2)y_3}{y_2} &\frac{1}{y_3} & (1 + (1 + y_1)y_2)y_3 \\
    y(2,m) & (1 + y_1)y_2 & \Delta &\frac{1 + y_3}{y_2y_3} & (1 + y_1)y_2 & \Delta \\
    y(3,m) & (1 + (1 + y_1)y_2)y_3 & \frac{1 + (1 + y_2)y_3}{y_1y_2y_3} & y_1 & \frac{1 + (1 + y_1)y_2}{y_1} & \frac{1 + (1 + y_2)y_3}{y_2}\\
    \hline
\end{array}
\]
$\Delta = \frac{1 + y_2 + (1 + (2 + y_1)y_2 + (1 + y_1)y_2^2)y_3}{y_1y_2}$.
    \caption{Expression of $y(i,m)'s$ as Laurent polynomials in the initial \bfY-cluster in type $A_3$.}
    \label{fig:yim-laurent}
\end{figure}
\endgroup

According to \eqref{eq:uc-to-Yfrieze}, the unordered \bfY-cluster $\{y_1, y_2, y_3\}$ gives rise to a well-defined arithmetic \bfY-frieze pattern which we denote by $k$. The first few values of $k$, obtained by setting $y_1 = y_2 = y_3 = 1$ in the Laurent expansion of each $y(i,m)$, are displayed below (c.f. \notaref{:infinite-grid} for the layout).

\begin{equation}
    \begin{tikzcd}[column sep=0em, row sep=0.25em]
&& &1&&3&&3&&1&&3&&3&&\\
	&&\cdots& &2&&8&&2&&2&&8&&2&&\cdots\\
    &&&& &3&&3&&1&&3&&3&&1&&
\end{tikzcd}\tag*{\text{$\diamond$}}
\end{equation}
\eex

\section{The \bfY-cluster structure of \bfY-frieze patterns}\label{s:clusteralg}

\subsection{Friezes of positive spaces}\label{ss:cluster-alg}
Let $B$ be any $r \times r$ mutation matrix, not necessarily $B = B_\sA$ and consider the Fock-Goncharov ensemble 
\[
\AA_B = (\FF(\AA_B), \; \{{\bf x}_t\}_{t \in \TT_r}) \hs \mbox{and} \hs \YY_B = (\FF(\YY_B), \; \{{\bf y}_t\}_{t \in \TT_r}),
\]
with common matrix pattern $\{B_t\}_{t \in \TT_r}$. Recall from \cite{GHKK,CdsgL:ringel-conj} that an element $x \in \calU(\AA_B)$ is called a {\it global monomial on} $\AA_B$ if there exist $t \in \TT_r$ and ${\bf m} \in \ZZ^r$ such that $x = {\bf x}_t^{\bf m}$. Similarly, an element $y \in \calU(\YY_B)$ is called a {\it global monomial on} $\YY_B$ if there exist $t \in \TT_r$ and ${\bf m} \in \ZZ^r$ such that $y = {\bf y}_t^{\bf m}$. It is well-known that the cluster monomials of $\AA_B$ are precisely all global monomials on $\AA_B$.

\bde{:cluster-alg}\cite[Definition 0.1]{GHKK}
Let $V = \AA_B$ or $\YY_B$.

1) The {\it cluster algebra $\calA(V)$ of $V$} is the sub-algebra of $\calU(V)$ generated by global monomials on $V$. 

2) A {\it frieze of $\calA(V)$} is a ring homomorphism 
\[
    \calA(V) \longrightarrow \ZZ
\]
taking positive integer values at every global monomial on $V$. \hfill $\diamond$
\ede

When $V = \AA_B$, this definition agrees with \cite[Definition 3.1]{GHL:friezes} and \cite[Definition 3.2]{GM:finite}. When $V = \YY_B$, this definition appears to be new. Examples of friezes of $\calA(\YY_B)$ when $B$ is of size 2 will be given in \S \ref{s:rank2}. For now, let us illustrate how friezes of $\calA(\AA_B)$ and friezes of $\calA(\YY_B)$ for the same $B$ can exhibit very different behaviours.

\bex{:Markov-Y}
Let $B = \begin{pmatrix}
    0 & 2 & -2 \\
    -2 & 0 & 2 \\
    2 & -2 & 0
\end{pmatrix}$ be the mutation matrix associated to the Markov quiver. On the one hand, there are infinitely many friezes of $\calA(\AA_B)$. On the other hand, there is only {\it one} frieze of $\calA(\YY_B)$. Indeed, one can check that the matrix pattern $\{B_t\}_{t \in \TT_3}$ obtained by setting $B_{t_0} = B$ contains exactly two matrices, $B$ and $-B$. Using \cite[Lemma 3.3]{CdsgL:ringel-conj}, we see that $\calA(\YY_B) \cong \ZZ[T^{\pm 1}]$. The single frieze of $\calA(\YY_B)$ is then given by 
\begin{equation}
    \ZZ[T^{\pm 1}] \longrightarrow \ZZ, \qquad T \longmapsto 1. \tag*{\text{$\diamond$}}
\end{equation}
\eex

Whilst $\calA(\AA_B)$ has been a central object in cluster theory since the theory's inception in the early 2000s, relatively little is known about $\calA(\YY_B)$. For this reason, we conclude this subsection by defining a sub-algebra of $\calA(V)$ which we believe ``approximates" $\calA(V)$ well. Whilst this sub-algebra will not play an important role in this article, we believe it may prove useful in the future.

Denote by $\mathfrak{X}(\AA_B)$ the set of cluster variables of $\AA_B$, and by $\mathfrak{X}(\YY_B)$ the set of \bfY-variables of $\YY_B$. Let $V = \AA_B$ or $\YY_B$. Recall \cite[\S 3.5]{CdsgL:ringel-conj} that elements in $\mathfrak{X}(V) \cap \calU(V)$ are called {\it global variables on $V$}. By the Laurent phenomenon, the cluster variables of $\AA_B$ are precisely {\it all} global variables on $\AA_B$. On the other hand, not every \bfY-variable of $\YY_B$ is a global variable on $\YY_B$.

\bde{:variable-alg}
Let $V = \AA_B$ or $\YY_B$. The {\it variable algebra $\calA^{\rm var}(V)$ of $V$} is the sub-algebra of $\calU(V)$ generated by the global variables on $V$. \hfill $\diamond$
\ede
For all mutation matrices, it is clear that $\calA^{\rm var}(\AA_B) = \calA(\AA_B) \subset \calU(\AA_B)$ and that
\begin{equation}\label{eq:var-alg-in-clus-alg}
   \calA^{\rm var}(\YY_B) \subset \calA(\YY_B) \subset \calU(\YY_B).
\end{equation}
In the setup of \exref{:Markov-Y}, one can show that the variable algebra is $\ZZ$. In other words, there exist mutation matrices $B$ such that the first inclusion in \eqref{eq:var-alg-in-clus-alg} is strict.

\subsection{Friezes and \bfY-frieze patterns}

We now relate the algebraic notions of friezes defined in \S \ref{ss:cluster-alg} to the elementary notions  of frieze patterns and \bfY-frieze patterns of \S \ref{s:frieze-patterns}. Let us begin once again by fixing a symmetrisable generalised Cartan matrix $A = (a_{i,j})$ and setting $B = B_\sA$ as in \eqref{eq:B_A}. Recall the definition of $x(i,m)$ and $y(i,m)$ from \S \ref{ss:unitary}.

\bde{:belt-alg}
 We continue with the notation above.

1) The {\it belt algebra $\calA^\flat(\AA_B)$} is the sub-algebra of $\calU(\AA_B)$ generated by $\{x(i,m)\}_{(i,m) \in [1,r] \times \ZZ}$. 

2) The {\it \bfY-belt algebra} $\calA^\flat(\YY_B)$ is the sub-algebra of $\calU(\YY_B)$ generated by $\{y(i,m)\}_{(i,m) \in [1,r] \times \ZZ}$. 
\ede

By definition, $\calA^\flat(\AA_B)$ is a sub-algebra of $\calA(\AA_B)$. The following proposition shows that the two algebras coincide.
\bpr{:cluster-alg-generator}
    The set $\{x(i,0), x(i,1) : i \in [1,r]\}$ generates $\calA(\AA_B)$ as a ring. In particular, $\calA^\flat(\AA_B) = \calA(\AA_B)$.
\epr
\begin{proof}
    Let $R = \ZZ[z_1, \ldots , z_{2r}]$ be the polynomial ring in $2r$ variables. The results in \cite[\S 6.3]{GHL:friezes} imply that there exists a {\it surjective} ring homomorphism 
    \[
        \psi^\BFZ: R \longrightarrow  \calA(\AA_B),
    \]
    such that $\psi^\BFZ(z_i) = x(i,0)$ and $\psi^\BFZ(z_{r+i}) = x(i,1)$ for all $i \in [1,r]$. The claim then follows.
\end{proof}

Suppose that $g$ is a frieze of $\calA(\AA_B)$. Setting 
\[
    f(i,m) = g(x(i,m)), \qquad \forall (i,m) \in [1,r] \times \ZZ,
\]
it is easy to see that $f \in {\rm Frieze}(A, \ZZ_{>0})$ (see \cite[\S 5.3]{GHL:friezes} for more). In this case we say that {\it f is obtained from $g$ by evaluation}. \prref{:cluster-alg-generator} and \cite[Proposition 5.11]{GHL:friezes} guarantee that the converse holds: {\it every $f \in {\rm Frieze}(A, \ZZ_{>0})$ is obtained from a frieze of $\calA(\AA_B)$ by evaluation}.

Turning now to the \bfY-spaces, recall that we have the following chain of inclusions:
\begin{equation}\label{eq:Y-alg-chain}
     \calA^\flat(\YY_B) \subset \calA^{\rm var}(\YY_B) \subset \calA(\YY_B) \subset \calU(\YY_B).
\end{equation}
When $B$ has size 2, we show in \leref{:Y-alg-chain-equal}, \coref{:cluster-equal-upper} that the inclusions in \eqref{eq:Y-alg-chain} are all equalities. When $B$ is of finite type, it was shown in \cite[Theorem B]{CdsgL:ringel-conj} the set $\{y(i,m) : (i,m) \in [1,r] \times \ZZ\}$ contains all global variables on $\YY_B$, and thus $\calA^\flat(\YY_B) =\calA^{\rm var}(\YY_B)$ in this case. It would be interesting to investigate the ``gaps" between the rings in \eqref{eq:Y-alg-chain}, if they exist. 

\bde{:Y-frieze}
By a {\it frieze of $\calA^\flat(\YY_B)$} we mean a ring homomorphism 
\[
    h: \calA^\flat(\YY_B) \longrightarrow \ZZ
\]
taking positive integer values on the set $\{y(i,m) : (i,m) \in [1,r] \times \ZZ\}$. \hfill $\diamond$
\ede

Suppose that $h$ is a frieze of $\calA^\flat(\YY_B)$. Setting 
\[
    k(i,m) = h(y(i,m)), \qquad \forall (i,m) \in [1,r] \times \ZZ
\]
and applying $h$ to both sides of \eqref{eq:y(m,i)-relations}, we see that $k \in {\rm YFrieze}(A, \ZZ_{>0})$. In this case, we say that $k$ is {\it obtained from $h$ by evaluation}. The following proposition says that this procedure establishes a bijection between \bfY-frieze patterns associated to $A$ and friezes of $\calA^\flat(\YY_B)$.

\bpr{:Y-frieze-pat-to-homo}
   Every arithmetic \bfY-frieze pattern associated to $A$ is obtained from a (necessarily unique) frieze of $\calA^\flat(\YY_B)$ by evaluation.
\epr
\begin{proof}
Uniqueness follows from the fact that the set $\{y(i,m) : (i,m) \in [1,r] \times \ZZ\}$ generates $\calA^\flat(\YY_B)$. It remains to show that every arithmetic \bfY-frieze pattern associated to $A$ arises in this way. Suppose that $f \in {\rm YFrieze}(A, \ZZ_{>0})$. Denote by ${\bf y} = (y_1, \ldots, y_r)$ the \bfY-cluster of $\YY_B$ at $t_0 = t(1,0)$. By a direct computation (see e.g. \cite[proof of Theorem 5.2]{CdsgL:ringel-conj}), we know that 
\[
    y_i = y(i,0) \prod_{j=1}^{i-1} (1+y(j,0))^{a_{j,i}}, \qquad i \in [1,r].
\]
Define 
\[
    h: \ZZ[y_1^{\pm 1}, \ldots , y_r^{\pm 1}] \longrightarrow \QQ, \qquad y_i \longmapsto f(i,0)\prod_{j=1}^{i-1} (1+f(j,0))^{a_{j,i}}.
\]
By \leref{:yim-universal-Laurent}, $h$ restricts to a ring homomorphism, also denoted by $h$, from $\calA^\flat(\YY_B)$ to $\QQ$ which takes positive values on every $y(i,m)$. It is not hard to check that
\[  
    h(y(i,0)) = f(i,0), \qquad i \in [1,r].
\]
By applying $h$ to both sides of \eqref{eq:y(m,i)-relations} and reasoning inductively with respect to the total order \eqref{eq:tot-order-grid}, we see that $  h(y(i,m)) = f(i,m)$ for all $(i,m)$.
Since $f$ is $\ZZ_{>0}$-valued by assumption, $h(\calA^\flat(\YY_B) ) \subset \ZZ$, i.e. $h$ is a frieze of $\calA^\flat(\YY_B)$ and we are done. 
\end{proof}

We conclude this section by establishing the connection between $p_\AS$ given in \S \ref{ss:friezepat} and the (pullback of the) ensemble map $p_{\!{}_B}^*$ defined in \S \ref{ss:patterns}. To do so, note that the semi-field homomorphism $p_{\!{}_B}^*: \FF_{>0}(\YY_B) \to \FF_{>0}(\AA_B)$ defined in \eqref{eq:p-map} uniquely determines a ring homomorphism $\calU(\YY_B) \to \calU(\AA_B)$, which we also denote by $p_{\!{}_B}^*$. It is straightforward to show that 
\[
    p_{\!{}_B}^*(y(i,m)) = {\bf x}_{t(i,m)}^{B_{t(i,m)}e_i} \stackrel{\rm denote}{=} {\hat y}(i,m), \qquad \forall (i,m) \in [1,r] \times \ZZ
\]
which by \eqref{eq:prop-Btim} and \eqref{eq:xt(i,m)} is a monomial in the variables $\{x(i,m) : (i,m) \in [1,r] \times \ZZ\}$. Thus,
$p_{\!{}_B}^*$ restricts to a ring homomorphism 
\[
\widetilde{p}_{\!{}_B}^{\,*} \stackrel{\rm denote}{=} p_{\!{}_B}^*\vert_{\calA^\flat(\YY_B)} \;: \; \calA^\flat(\YY_B) \longrightarrow \calA(\AA_B).
\]

\bco{:pas-and-ensemb-map}
Let $S = \ZZ_{>0}$ and $f \in {\rm Frieze}(A,S)$. If $f$ is obtained from the frieze $g: \calA(\AA_B) \to \ZZ$ by evaluation, then $p_\AS(f)$ is obtained from $g \circ \widetilde{p}_{\!{}_B}^{\,*}: \calA^\flat(\YY_B) \to \ZZ$ via evaluation.
\eco

\section{\bfY-frieze patterns in rank 2}\label{s:rank2}
\subsection{Background on generalised cluster algebras}\label{ss:gen-clus-alg}
In this subsection, we recall some facts about generalised cluster algebras of rank 2. For a more detailed discussion, see \cite{rupel:greedy-gen-clust-alg}. 

Fix $a,b$ two positive integers. Let $n$ be a positive integer, $t$ a formal variable and set $P_n(t) = (1+t)^n$. It is clear that $P_n$ is {\it monic} and that 
$t^n P(t^{-1}) = P(t)$, i.e. $P_n$ is {\it palindromic}. Let $x_1,x_2$ be two formal commuting variables, and define $\{x_k\}_{k \in \ZZ}$ recursively by the relations
\begin{equation}
    x_k x_{k+2} = \begin{cases}
        P_c(x_{k+1}), \qquad&k\text{ odd,} \\
        P_b(x_{k+1}), \qquad&k \text{ even.}
    \end{cases}
\end{equation}
Consider
\[
    \calU(P_b,P_c) = \bigcap_{k\in \ZZ} \ZZ[x_{k}^{\pm 1}, x_{k+1}^{\pm 1}] \quad \text{and} \quad \calU^+(P_b,P_c) = \bigcap_{k\in \ZZ} \ZZ_{\geq 0}[x_{k}^{\pm 1}, x_{k+1}^{\pm 1}]
\]
and denote by $\calA(P_b,P_c)$ the sub-algebra of $\QQ(x_1,x_2)$ generated by $\{x_k\}_{k \in \ZZ}$. The algebra $\calA(P_b,P_c)$ is a {\it generalised cluster algebra} in the sense of \cite{rupel:greedy-gen-clust-alg}. Each $x_k$ is called a {\it generalised cluster variable} and each pair $(x_k, x_{k+1})$ a {\it generalised cluster} of $\calA(P_b,P_c)$. Monomials in variables in a generalised cluster are called {\it generalised cluster monomials} of $\calA(P_b,P_c)$. 

\bex{:gca-finite-type}
1) Let $b=c = 1$. Then $\calA(P_b, P_c)$ is the (ordinary) cluster algebra of type $A_2$. 

2) Let $b = 2$, $c = 1$. Then $x_{k} = x_{k+6}$ for all $k \in \ZZ$, and
\[
\begin{aligned}
    x_3 &= \frac{x_2+1}{x_1}, &     x_5&= \frac{x_1^2+2x_1+x_2+1}{x_1x_2},\\
    x_4 &= \frac{(1+x_1+x_2)^2}{x_1^2x_2}, & x_6 &= \frac{(1+x_1)^2}{x_2}.
\end{aligned}
\]
3) Let $b =3$, $c = 1$. Then $x_{k} = x_{k+8}$ for all $k \in \ZZ$, and 
\[
\begin{aligned}
    x_3 &= \frac{1+x_2}{x_1}, &  x_5 &=\frac{(1+x_1)^3+3x_1x_2+x_2^2+2x_2}{x_1^2x_2} ,&   x_7 &= \frac{(1+x_1)^3+x_2}{x_1x_2}, \\
     x_4 &= \frac{(1+x_1+x_2)^3}{x_1^3x_2}, &  x_6 &= \frac{(x_1^2+2x_1+x_2+1)^3}{x_1^3x_2^2},&   x_8 &= \frac{(1+x_1)^3}{x_2}. 
\end{aligned}
\]
\eex

The following is a summary of results we will need.

\ble{:gca-properties}\cite[Theorems 1.1, 2.4, 2.5]{rupel:greedy-gen-clust-alg}
The following statements hold. 

1) $\calA(P_b,P_c) = \calU(P_b,P_c)$.

2) For each $k \in \ZZ$, $\{x_k, x_{k+1}, x_{k+2}, x_{k+3}\}$ is a generating set of $\calA(P_b,P_c)$ as a ring.

3) For each $k \in \ZZ$, $x_k \in \calU^+(P_b,P_c)$. 
\ele

We now introduce friezes of $\calA(P_b,P_c)$. 

\bde{:gca-frieze}
    A {\it frieze} of $\calA(P_b,P_c)$ is a ring homomorphism  $       \calA(P_b,P_c) \longrightarrow \ZZ$ 
    which takes positive integer values at all generalised cluster variables of $\calA(P_b,P_c)$.  \hfill $\diamond$
\ede

The following lemma will be used in the proof of \thref{:Yfrieze-count-rk2}.
\ble{:gencluster-to-frieze}
Given a generalised cluster ${\bf x} = (x_1, x_2)$ of $\calA(P_b,P_c)$, the ring homomorphism 
\[
    \psi_{\bf x}: \ZZ[{\bf x}^{\pm 1}] \longrightarrow \ZZ, \quad \psi_{\bf x}(x_i) = 1, \; i = 1,2,
\]
restricts to a frieze of $\calA(P_b,P_c)$. Moreover, distinct (as unordered sets) generalised clusters give rise to distinct friezes of $\calA(P_b,P_c)$. 
\ele
\begin{proof}
    The restriction of $\psi_{\bf x}$ to $\calA(P_b,P_c)$ is obviously a frieze of $\calA(P_b,P_c)$, thanks to 3) in \leref{:gca-properties}. Now suppose that ${\bf x}= (x_1,x_2)$ and ${\bf z} = (z_1,z_2)$ are two distinct generalised clusters such that $\psi_{\bf x} = \psi_{\bf z}$. Without loss of generality, we may assume that $z_1 \not\in \{x_1,x_2\}$. On the one hand, $\psi_{\bf x}(z_1) = 1$, and on the other $z_1 \in \ZZ_{\geq 0}[x_1^{\pm 1}, x_2^{\pm 1}]$ by \leref{:gca-properties}. This implies that $z_1$ is a Laurent monomial of $x_1$ and $x_2$, say 
    \[
        z_1 = x_1^{d_1} x_2^{d_2}.
    \]
    If  $d_1,d_2 \geq 0$, i.e. $x_1^{d_1} x_2^{d_2}$ is a generalised cluster monomial, then $z_1 - x_1^{d_1} x_2^{d_2} = 0$, which contradicts the linear independence of generalised cluster monomials proved in \cite[Theorem 1.1]{rupel:greedy-gen-clust-alg}. Thus we may assume that $d_1< 0$ or $d_2 <0$. Without loss of generality, say $d_1 < 0$. Let ${\bf x}'$ denote the unique generalised cluster of $\calA(P_b,P_c)$ containing $x_2$ but not $x_1$. By a straightforward argument, one can show that $z_1$ is not a Laurent polynomial in the generalised cluster variables of ${\bf x}'$, contradicting \leref{:gca-properties}. 
\end{proof}

\subsection{Arithmetic \bfY-frieze patterns in rank 2}
Let $b,c$ be as in \S \ref{ss:gen-clus-alg}, and consider the $2 \times 2$ symmetrisable generalised Cartan matrix 
\begin{equation}\label{eq:cartan-rank2}
    A = \begin{pmatrix}
        2 & -b \\ -c & 2
    \end{pmatrix}.
\end{equation}
Note that since $b, c \neq 0$ by assumption, the matrix $A$ is indecomposable in the sense of \S \ref{ss:notation}. We set $B = B_\sA$ as in \eqref{eq:B_A}, and consider the positive space
\[
    \YY_B = (\FF(\YY_B), \{{\bf y}_t = (y_{t;1}, y_{t;2})\}_{t \in \TT_2}).
\]
Setting $t_0 = t(1,0)$ as in \S \ref{ss:unitary}, one sees that every $t \in \TT_2$ is of the form $t(i,m)$ for a unique pair $(i,m) \in [1,2]\times \ZZ$.

Recall that we have defined the rings $\calA^\flat(\YY_B)$, $\calA^{\rm var}(\YY_B)$, $\calA(\YY_B)$ and $\calU(\YY_B)$ associated to $\YY_B$ in \S \ref{ss:cluster-alg} and \deref{:belt-alg}. 

\ble{:Y-alg-chain-equal}
With the notation as above, we have 
\[
    \calA^\flat(\YY_B) = \calA^{\rm var}(\YY_B) = \calA(\YY_B).
\]
\ele
\begin{proof}
    It follows easily from \cite[Lemma 3.3]{CdsgL:ringel-conj} that the set of global variables on $\YY_B$ is precisely $\{y(i,m): (i,m) \in [1,r] \times \ZZ\}$, and thus $ \calA^\flat(\YY_B) = \calA^{\rm var}(\YY_B)$, by definition. Moreover, the set of global monomials on $\YY_B$ is readily seen to be
    \begin{equation}\label{eq:global-rank2}
        \{{\bf y}_{t(1,m)}^{\bf a}, {\bf y}_{t(2,m)}^{-{\bf a}} : {\bf a} = (a_1, -a_2)^T, \; a_1,a_2 \in \ZZ_{\geq 0} \}.
    \end{equation}
    Using \eqref{eq:Y-mut} and the definition of $y(i,m)$ given in \S \ref{ss:unitary}, we have
    \begin{equation}\label{eq:Ycluster-rank-2}
        {\bf y}_{t(1,m)} = (y(1,m), y(2,m-1)^{-1}) \quad\text{and}\quad {\bf y}_{t(2,m)} = (y(1,m)^{-1}, y(2,m)), \qquad \forall m \in \ZZ.
    \end{equation}
    By combining \eqref{eq:global-rank2} and \eqref{eq:Ycluster-rank-2}, every global monomial on $\YY_B$ is a monomial in $\{y(i,m): (i,m) \in [1,2] \times \ZZ\}$, and therefore $\calA(\YY_B) \subset \calA^\flat(\YY_B)$. The claim then follows from \eqref{eq:var-alg-in-clus-alg}.
\end{proof}

\bpr{:Y-alg-is-gca}
The field isomorphism 
\begin{equation}\label{eq:Y-to-gca}
    \FF(\YY_B) \longrightarrow \QQ(x_1,x_2), \quad y(1,0) \mapsto x_1, \; y(2,0) \mapsto x_2,
\end{equation}
restricts to an isomorphism of rings 
\[
    \phi: \calA(\YY_B) \longrightarrow\calA(P_b, P_c).
\]
Moreover, $\phi$ maps global monomials on $\YY_B$ to generalised cluster monomials of $\calA(P_b, P_c)$, and global variables on $\YY_B$ to generalised cluster variables of $\calA(P_b, P_c)$. In particular, friezes of $\calA^\flat(\YY_B) $ are in bijection with friezes of $\calA(P_b, P_c)$ via pull-back by $\phi$. 
\epr
\begin{proof}
    By \leref{:Y-alg-chain-equal}, it is enough to determine the image of the $y(i,m)$'s. We show by induction on the total order \eqref{eq:tot-order-grid} that 
    \begin{equation}\label{eq:y-to-var}
        \phi(y(i,m)) = x_{2m+i}
    \end{equation} 
    for all $m \in \ZZ$ and $i = 1,2$. 
    By definition, $\phi(y(i,0)) = x_i$ for $i = 1,2$. Now fix $m \geq 0$, and assume that \eqref{eq:y-to-var} holds for $i = 1,2$ and $0 \leq k \leq m$. Then, 
    \[
        \phi(y(1,m+1)) = \phi\left(\frac{(1+y(2,m))^c}{y(1,m)}\right) = \frac{(1+x_{2m+2})^c}{x_{2m+1}} = x_{2(m+1)+1}.
    \]
    Similarly, $\phi(y(2,m+1)) = x_{2(m+1)+2}$. Thus \eqref{eq:y-to-var} holds for all $m \geq 0$. By a similar argument, one can show that \eqref{eq:y-to-var} holds for all $m \in \ZZ$. In particular, $\phi$ maps $\calA(\YY_B)$ bijectively onto $\calA(P_b,P_c)$. The remaining claims follow by combining \eqref{eq:y-to-var} with \eqref{eq:global-rank2}.
\end{proof}

\bco{:cluster-equal-upper}
With the notation as above, $\calA(\YY_B) = \calU(\YY_B)$.
\eco
\begin{proof}
    Using \eqref{eq:Ycluster-rank-2}, we see that 
    \[
        \calU(\YY_B) = \bigcap_{m \in \ZZ} \left(\ZZ[y(1,m)^{\pm 1}, y(2,m-1)^{\pm 1}] \cap \ZZ[y(1,m)^{\pm 1}, y(2,m)^{\pm 1}]\right).
    \]
    Moreover, \eqref{eq:y-to-var} implies that the image of $\calU(\YY_B)$ under \eqref{eq:Y-to-gca} is precisely $\calU(P_b,P_c) = \calA(P_b,P_c)$ (this last equality is 1) of \leref{:gca-properties}). The result then follows from \prref{:Y-alg-is-gca}. 
\end{proof}

\begin{figure}
    \centering
    \includegraphics[width=1\textwidth]{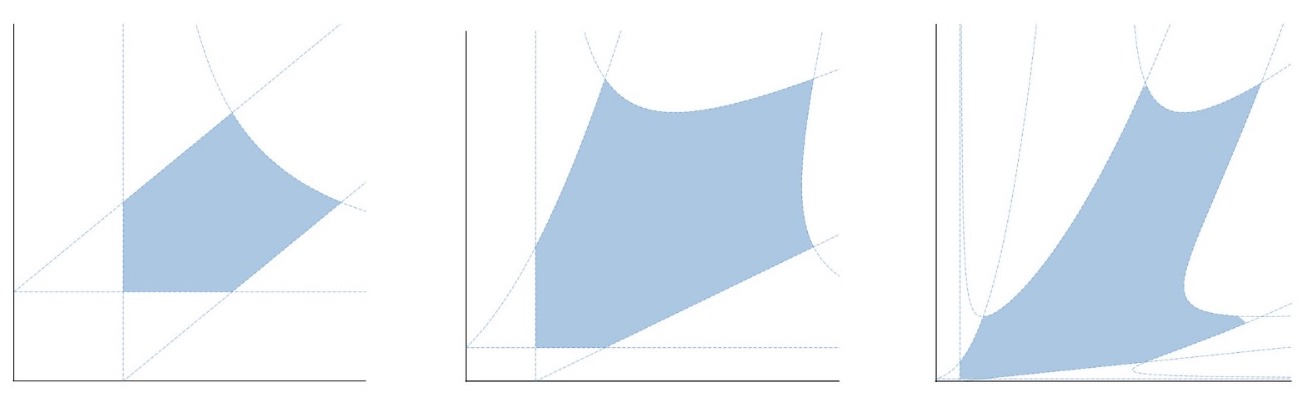}
    \caption{The superunitary region (embedded in $\RR^2$) in types $A_2, C_2$ and $G_2$.}
    \label{fig:superunitary-finite}
\end{figure}

For convenience, we write $\calA = \calA(P_b, P_c)$. In analogy with the discussion in \cite[\S 1.2]{GM:finite}, the {\it superunitary region} of $\calA$, denoted $\calA(\RR_{\geq 1})$, is the set of all ring homomorphisms 
\[
    \calA \longrightarrow \RR
\]
which take values in $\RR_{\geq 1}$ for all generalised cluster variables of $\calA$. It is clear that 
\[
    \calA(\RR_{\geq 1}) \longrightarrow (\RR_{\geq 1})^2, f \longmapsto (f(x_1), f(x_2)),
\]
is an embedding of $\calA(\RR_{\geq 1})$ into $(\RR_{\geq 1})^2 \subset (\RR_{>0})^2$. Identifying $\calA(\RR_{\geq 1})$ with its image in $(\RR_{>0})^2$, the superunitary region can described explicitly by a (potentially infinite set) of polynomial inequalities, one for each generalised cluster variable of $\calA$ (expressed as a Laurent polynomial in $x_1$ and $x_2$). Using \exref{:gca-finite-type}, we graph the superunitary region (embedded in $\RR^2$) when $1 \leq bc \leq 3$ in Figure \ref{fig:superunitary-finite}.

\bth{:Yfrieze-count-rk2}
Let $A$ be any $2 \times 2$ indecomposable generalised symmetrisable Cartan matrix. There are finitely many arithmetic \bfY-frieze patterns associated to $A$ if and only if $A$ is of finite type. When $A$ is of finite type, the count of arithmetic \bfY-frieze patterns is given in the following table.
\begingroup
\renewcommand{\arraystretch}{1.5}
\[
\begin{array}{|c|c|c|c|}
\hline 
   A   & A_2 & C_2 & G_2 \\ \hline
   \text{\# of arithmetic \bfY-frieze patterns}  & 5 & 10 & 21 \\
   \hline
\end{array}
\]
\endgroup

\eth
\begin{proof}
    Suppose that $A$ is indecomposable of finite type, i.e. A is of the form \eqref{eq:cartan-rank2} with $1 \leq bc \leq 3$. By computing the generalised cluster variables of $\calA = \calA(P_b,P_c)$ explicitly as in \exref{:gca-finite-type}, one sees that the superunitary region $\calA(\RR_{\geq 1})$ of $\calA$, embedded in $\RR^2$, is a compact set; see Figure \ref{fig:superunitary-finite}. This implies that the set of friezes of $\calA$, which is a closed and discrete subset of $\calA(\RR_{\geq 1})$, is finite. By \prref{:Y-alg-is-gca} and \prref{:Y-frieze-pat-to-homo}, the set of arithmetic \bfY-frieze patterns associated to $A$ is finite. The number of arithmetic \bfY-frieze patterns can then be calculated using a computer program. 

    Now suppose that $A$ is of infinite type, i.e. $bc \geq 4$. It is well-known \cite{FZ:ClusterII, CS:gen-cluster-alg} that $\calA(P_b,P_c)$ has infinitely many generalised clusters, and by \leref{:gencluster-to-frieze} there are infinitely many friezes of $\calA(P_b,P_c)$. The claim that there are infinitely many \bfY-frieze patterns assocaited to $A$ then follows from the bijections in \prref{:Y-alg-is-gca} and \prref{:Y-frieze-pat-to-homo}.  
\end{proof}

\section{{\bf Y}-frieze patterns in finite type}\label{s:enum-problem}

\subsection{Symmetries of {\bf Y}-frieze patterns}\label{ss:glide-sym}
Fix $A = (a_{i,j})$ an $r \times r$ indecomposable Cartan matrix of finite type, and let ${\bf h}$ denote the Coxeter number associated to $A$. We recall some notation introduced in \cite[\S 7.5]{GHL:friezes}. 

1) When the Cartan matrix $A$ is of type $A_1, B,C,D_{2n}, n \geq 2,  E_7,E_8, F_4$ or $G_2$, set  
\[
   i = i^* \quad \text{ and } \quad  m_i = \frac{{\bf h}}{2} , \qquad  i \in [1, r].
 \]

2) When $A$ is of type $A_n, D_{2n+1}$, $n \geq 2$, or $E_6$, define the quiver $Q$ whose vertex set is $[1, r]$ and where, for each pair $i \neq j$ there is one arrow from $i$ to $j$ if and only if 
$i < j$ and $a_{i, j} =-1$. The map  $i \mapsto i^*$ is the permutation of the vertex set $[1,r]$ induced by the unique graph automorphism of the Dynkin diagram underlying $Q$. For
 $i \in [1,r]$, let 
 \[
 \pi(i, i^*): \; \;i = i_0 \mathdash \;i_1 \mathdash 
 \;i_2 \mathdash  \cdots \mathdash \; i_{l-1}  \mathdash \;i_l = i^*
 \]
 be the {\it unique path in the Dynkin diagram connecting $i$ and $i^*$}. If $i = i^*$, we set $a_i = b_i = 0$. If $i \neq i^*$, 
  let $a_i$ be the 
number of arrows in $Q$, whose underlying edges lie in $\pi(i, i^*)$, that are pointing towards $i^*$, and let $b_i$ be the number of arrows in $Q$, whose underlying edges lie in $\pi(i, i^*)$, that are pointing towards $i$. For $i \in [1,r]$, set
\[
m_i = \frac{1}{2}({\bf h} + a_i - b_i).
\]

Define the bijection (see \cite[\S 7.5 and Remark 7.25]{GHL:friezes})
\[
{\bf F}_\sA: [1,r] \times \ZZ \longrightarrow [1,r] \times \ZZ, \qquad {\bf F}_\sA(i,m) = (i^*, m+ m_{i^*} + 1).
\]

 Recall from \reref{:frieze-pattern-history} and \reref{:y-frieze-history} the definitions of the generic frieze pattern and the generic \bfY-frieze pattern associated to $A$.
\ble{:generic-glide}\cite[Proposition 6.2]{CdsgL:ringel-conj}
    The generic frieze pattern and the generic \bfY-frieze pattern associated to $A$ are ${\bf F}_\sA$-invariant.
\ele

Let $S$ be a semi-ring, and consider the statement: 
\begin{equation}\label{eq:embedding-condition}
    S \text{ is the sub-semiring of some semifield}.
\end{equation} 
All semi-rings and semi-fields encountered in \exref{:semi-rings} and \exref{:semi-fields} satisfy \eqref{eq:embedding-condition}. 
\bth{:y-frieze-glide-symmetry}
    For any semi-ring $S$ satisfying \eqref{eq:embedding-condition}, $S$-valued frieze patterns and $S$-valued \bfY-frieze patterns associated to $A$ are ${\bf F}_\sA$-invariant.
\eth
\begin{proof}
We prove the claim for \bfY-frieze patterns. To begin, assume that $S$ is a semi-field. Recall from \leref{:enum-Yfrieze-semifield} that any $S$-valued \bfY-frieze pattern is of the form $k_{\bf s}$ for some ${\bf s} = (s_1, \ldots , s_r) \in S^r$. Let $\phi_s: \QQ_{>0}(Y_1, \ldots , Y_r)$ denote the unique (see \cite[Lemma 2.1.6]{BFZ:param-canon-basis}) semi-field homomorphism given by 
\[
    \phi_{\bf s}(Y_i) = s_i, \qquad i \in [1,r],
\]
and let $G$ denote the generic \bfY-frieze pattern. For all $(i,m)$, we have 
\[
    k_{\bf s}(i,m) = \phi_{\bf s} \circ G(i,m).
\]
The claim then follows from \leref{:generic-glide}. Suppose now that $S$ is a semi-ring satisfying \eqref{eq:embedding-condition}, and let $\overline{S}$ be a semi-field such that $S \subset \overline{S}$, as semi-rings. Since every $S$-valued \bfY-frieze pattern is in particular an $\overline{S}$-valued \bfY-frieze pattern, the claim follows. 

The claim for $S$-valued frieze patterns is obtained {\it mutatis mutandis}, and we omit it. 
\end{proof}

\subsection{Enumerating {\bf Y}-frieze patterns}\label{ss:enumprob}

When $S$ is a semi-field, $S$-valued \bfY-frieze patterns associated to $A$ are in bijection with  elements in $S^r$ (c.f. \leref{:enum-Yfrieze-semifield}). Note in particular that this is independent of the choice of symmetrisable generalised Cartan matrix. 

On the other hand, one does not expect such a parametrisation of $S$-valued \bfY-frieze patterns, when $S$ is a semi-ring {\it which is not a semi-field} (as pointed out at the end of \S \ref{ss:Yfriezepat}). In this case, determining the number of $S$-valued \bfY-frieze patterns is in general rather delicate. In this section, we study this problem when $S$ is the tropical semi-ring (\prref{:friezes-enum}) or the semi-ring of positive integers (\thref{:inf-Yfrieze}).

\bpr{:friezes-enum}
Let $A$ be an indecomposable Cartan matrix of finite type. Then, there are no non-trivial $\ZZ^{\max}_{\geq 0}$-valued \bfY-frieze patterns associated to $A$. 
\epr
Before providing a proof, we mention the following interesting consequence. 
\bco{:pas-bij-tropical-sr}
Let $A$ be an indecomposable Cartan matrix of finite type, and let $S = \ZZ^{\max}_{\geq 0}$. Then $p_\AS$ is bijective.
\eco
The proof of \prref{:friezes-enum} is equivalent to the well-known statement that there are no non-trivial $\ZZ_{\geq 0}$-valued additive functions (c.f. \cite{CdsgL:ringel-conj} and \cite[\S 1]{Ringel:cluster-additive}) associated to $A$, when $A$ is of finite type. Unfortunately, we were unable to find a proof in the literature which covered the non-simply laced case; we provide a proof for completeness. To do so, we need to introduce some notation coming from Lie theory. First, let 
$G$ be the connected and simply connected 
complex simple Lie group with same Cartan-Killing type as $A$. Fix a maximal torus $T$ and a Borel subgroup 
$B$ containing $T$, and let
$\calP$ be the character lattice of $T$. Let $\{\alpha_1, \ldots, \alpha_r\} \subset \calP$ be the set of simple roots determined 
by the pair $(B, T)$, and let $\Phi_+\subset \calP$ be the corresponding set of positive roots.

Then, the set $\{\alpha_1^\vee, \ldots, \alpha_r^\vee\}$ of co-roots in the co-character lattice of $T$ are such that
\begin{equation}\label{eq:Cartan-basis-change}
    a_{i, j} = (\alpha_i^\vee, \alpha_j), \hs \forall\;\; i, j \in [1, r].
\end{equation}
Since $G$ is simply connected, $\{\alpha_1^\vee, \ldots, \alpha_r^\vee\}$ is a basis of the co-character lattice of $T$.
Let $\{\omega_1, \ldots, \omega_r\}$ be the basis of $\calP$ dual to $\{\alpha_1^\vee, \ldots, \alpha_r^\vee\}$, so that $\underline{\alpha} = \underline{\omega} A$,
where $\underline{\alpha} =  (\alpha_1, \ldots, \alpha_r)$ and $\underline{\omega} =  (\omega_1, \ldots, \omega_r)$. Let $W$ be the Weyl group generated by the reflections $s_i = s_{\alpha_i}$
for $i \in [1, r]$, and consider the Coxeter element $
c = s_1s_2 \cdots s_r$. Introduce $\underline{\beta} = (\beta_1, \ldots, \beta_r)$, where
\begin{equation}\label{eq:beta-i}
\beta_i = s_1 s_2 \cdots s_{i-1}(\alpha_i) \in \Phi_+, \hs i \in [1, r].
\end{equation}
Then $\{\beta_i: i \in [1, r]\} = \{\beta \in \Phi_+: c^{-1}\beta \in -\Phi_+\}$.
Write $C = (-L_\sA U_\sA^{-1})^T$, and
\begin{equation}\label{eq:L0U0}
L_\sA = \left(\begin{array}{cccc} 1 & 0  & \cdots & 0\\
a_{2,1} & 1 & \cdots & 0\\
\cdots & \cdots & \cdots& \cdots\\
\vspace{.08in}
a_{r,1} & a_{r,2} & \cdots &1\end{array}\right) \hs\mbox{and} \hs
U_\sA = \left(\begin{array}{cccc} 1 & a_{1,2}  & \cdots & a_{1,r}\\
0 & 1 & \cdots & a_{2,r}\\
\cdots & \cdots & \cdots& \cdots\\
\vspace{.08in}
0 & 0 & \cdots &1\end{array}\right).
\end{equation}
It is straightforward to check that \cite[Lemma 2.1]{YZ:Lcc}
\begin{equation}\label{eq:beta-omega}
\underline{\alpha} = \underline{\beta} \,U_\sA, \quad c\underline{\alpha} = - \underline{\beta} L_\sA \hs \mbox{and} \hs \underline{\beta} = (1-c) \underline{\omega}.
\end{equation}
Let $\calQ ={\rm Span}_\ZZ\{\alpha_1, \ldots, \alpha_r\} \subset \calP$ be the root lattice. It follows from \eqref{eq:beta-omega}
that $\underline{\beta}$ is a $\ZZ$-basis of $\calQ$ and that $1-c: \calP \to \calQ$ is a lattice isomorphism. It also follows from \eqref{eq:beta-omega} that
\begin{equation}\label{eq:c-omega}
c\, \underline{\beta} = -\underline{\beta}\, C^T.
\end{equation}

Since $c$ has order ${\bf h}$ and does not
have an eigenvalue equal to 1, we deduce that 
\begin{equation}\label{eq:charac-poly-Coxeter}
	I_r + C + \dots + C^{{\bf h} - 1} = 0.
\end{equation}

\begin{proof}[Proof of \prref{:friezes-enum}]
By definition, a $\ZZ_{\geq 0}$-valued \bfY-frieze pattern associated to $A$ is an assignment $f: [1,r] \times \ZZ \to \ZZ_{\geq 0}$ such that 
\begin{equation}\label{eq:tropical-Y-frieze}
    f(i,m) + f(i,m+1) = \sum_{j=i+1}^r (-a_{j,i}) f(j,m) + \sum_{j=1}^{i-1}(-a_{j,i})f(j,m+1),
\end{equation}
for each $(i,m) \in [1,r] \times \ZZ$. Clearly, the map (called the {\it trivial map}) 
\[
(i,m) \mapsto 0, \quad \forall (i,m) \in [1,r] \times \ZZ
\]
is a $\ZZ_{\geq 0}$-valued \bfY-frieze pattern associated to $A$. Now, suppose that $f$ is a $\ZZ_{\geq 0}$-valued \bfY-frieze pattern associated to $A$. Setting $
    f_m = (f(1,m), \ldots , f(r,m))^T \in (\ZZ_{\geq 0})^r$ and re-arranging \eqref{eq:tropical-Y-frieze}, we have
\[
    L_\sA^T f_m = -U_\sA^T f_{m+1}.
\]
Taking $(-U_\sA^T)^{-1}$ on both sides, we obtain $f_m = C^m f_0$ for all $m \in \ZZ$. Now, \eqref{eq:charac-poly-Coxeter} implies that
\[
    0 = (I_r+C + C^2 + \cdots + C^{{\bf h}-1})f_0 = f_0 + f_1 + f_2 + \cdots + f_{{\bf h}-1}.
\]
Since $f$ is $\ZZ_{\geq 0}$-valued, the equality holds if and only if all entries of $f_0, f_1, \ldots $ are zero, i.e. $f$ is the trivial map.
\end{proof}

We now turn to arithmetic \bfY-frieze patterns associated to $A$. The proof of the following theorem is a simple adaptation of that in \cite{muller:finite}.
\bth{:inf-Yfrieze}
    Let $A$ be an indecomposable Cartan matrix of finite type. The set of arithmetic $\bfY$-frieze patterns associated to $A$ is finite.
\eth
\begin{proof}
  By \thref{:y-frieze-glide-symmetry}, there exists a positive integer (depending on $A$) such that, for all $k \in {\rm YFrieze}(A,\ZZ_{>0})$ and all $(i,m) \in [1,r] \times \ZZ$, we have 
  \begin{equation}
      k(i,m) = k(i,m+p). 
  \end{equation}
  For each $i \in [1,r]$, set $K_i = \prod_{m=1}^p k(i,m)$. Since $k$ is $\ZZ_{>0}$-valued, $K_i \geq 1$. We now show that $K_i$ is bounded above. To do so, consider 
  \[
    K_i^2 = \prod_{m=1}^p k(i,m)^2 = \prod_{m=1}^p k(i,m) k(i,m+1).
  \]
  Substituting \eqref{eq:Yfrieze-rel} in each term in the product on the right-hand side, we obtain 
  \[
    K_i^2 = \prod_{m=1}^p \prod_{j=i+1}^r (1+k(j,m))^{-a_{j,i}} \prod_{j=1}^{i-1} (1+k(j,m+1))^{-a_{j,i}}. 
  \]
  Write 
  \[
    M_i = \prod_{m=1}^p \prod_{j=i+1}^r k(j,m)^{-a_{j,i}} \prod_{j=1}^{i-1} k(j,m+1)^{-a_{j,i}} = \prod_{j \neq i} K_j^{-a_{j,i}}.
  \]
  Set $b_i = \prod_{j \neq i}2^{-a_{j,i}}$. Using the identity $x < x+1 \leq 2x$ which holds for all $x \geq 1$, we deduce that 
  \begin{equation}\label{eq:bound}
      M_i < K_i^2 \leq b_i^p M_i.
  \end{equation}
  Set ${\bf K} = (K_1, \ldots , K_i, \ldots , K_r)$. By dividing out by $M_i$ in \eqref{eq:bound} for each $i \in [1,r]$, we obtain (using the notation in \S \ref{ss:notation}) that $ {\bf K}^A \leq (b_1^p, \ldots , b_r^p)$. It is clear from \eqref{eq:Cartan-basis-change} that $A$ is invertible and the entries of $A^{-1}$ are all positive rationals (see e.g. \cite[Table 2]{vinberg:liegroup}). In particular, we deduce that 
  \[
    {\bf K} \leq (b_1^p, \ldots , b_r^p)^{A^{-1}},
  \]
  i.e. each $K_i$ is bounded above. 
  Since $k$ is $\ZZ_{>0}$-valued, this implies that each $k(i,m)$ is bounded above, which concludes the proof.
\end{proof}

\bibliographystyle{alpha}
\bibliography{ref}

\begin{thebibliography}{dSGHL23}

\bibitem[ARS10]{ARS:friezes}
I.~Assem, C.~Reutenauer, and D.~Smith.
\newblock Friezes.
\newblock {\em Adv. Math.}, 225(6):3134--3165, 2010.

\bibitem[BFZ96]{BFZ:param-canon-basis}
A.~Berenstein, S.~Fomin, and A.~Zelevinsky.
\newblock Parametrizations of canonical bases and totally positive matrices.
\newblock {\em Adv. Math.}, 122(1):49--149, 1996.

\bibitem[CC73]{Con-Cox1}
J.H. Conway and H.S.M. Coxeter.
\newblock Triangulated polygons and frieze patterns.
\newblock {\em Math. Gaz.}, 57(400):87--94, 1973.

\bibitem[CdSGL23]{CdsgL:ringel-conj}
P.~Cao, A.~de~Saint~Germain, and J.-H. Lu.
\newblock {Tropical patterns from the Fock-Goncharov pairing and Ringel's conjectures}.
\newblock arXiv: 2311.17712, 2023.

\bibitem[Cox71]{Cox}
H.S.M. Coxeter.
\newblock Frieze patterns.
\newblock {\em Acta Arith.}, 18:297--310, 1971.

\bibitem[CS14]{CS:gen-cluster-alg}
L.~Chekhov and M.~Shapiro.
\newblock {Teichmüller spaces of Riemann surfaces with orbifold points of arbitrary order and cluster variables}.
\newblock {\em Int. Math. Res. Not. IMRN}, (10):2746–2772, 2014.

\bibitem[dSGHL23]{GHL:friezes}
A.~de~St.~Germain, M.~Huang, and J.H. Lu.
\newblock Friezes of cluster algebras of geometric type.
\newblock arXiv: 2309.00906, 2023.

\bibitem[FG09]{FG:ensembles}
V.~Fock and A.~Goncharov.
\newblock {Cluster ensembles, quantization and the dilogarithm}.
\newblock {\em Ann. Sci. de l’ENS}, 42(6):865--930, 2009.

\bibitem[FZ03]{FZ:ClusterII}
S.~Fomin and A.~Zelevinsky.
\newblock {Cluster algebras. II. Finite type classification}.
\newblock {\em Invent. Math}, 154(1):63--121, 2003.

\bibitem[FZ07]{FZ:ClusterIV}
S.~Fomin and A.~Zelevinsky.
\newblock {Cluster algebras. IV. Coefficients}.
\newblock {\em Compos. Math.}, 143(1):112--164, 2007.

\bibitem[GHKK18]{GHKK}
M.~Gross, P.~Hacking, S.~Keel, and M.~Kontsevich.
\newblock {Canonical bases for cluster algebras}.
\newblock {\em J. Amer. Math. Soc.}, 31(2):497--608, 2018.

\bibitem[GM22]{GM:finite}
E.~Gunawan and G.~Muller.
\newblock Superunitary regions of cluster algebras.
\newblock arXiv: 2208.14521v1, 2022.

\bibitem[Guo13]{Guo:tropical-frieze}
L.~Guo.
\newblock {On tropical friezes associated with Dynkin diagrams}.
\newblock {\em Int. Math. Res. Not. IMRN}, 2013(18):4243--4284, 2013.

\bibitem[Kac83]{Kac:inf-dim-Lie-alg}
V.G. Kac.
\newblock {\em Infinite-dimensional Lie algebras. An introduction.}, volume~44 of {\em Progress in Mathematics}.
\newblock Birkhäuser Boston, Inc., Boston, MA, 1983.

\bibitem[Kel10]{keller-knitting-algo}
B.~Keller.
\newblock {Cluster algebras, quiver representations and triangulated categories}.
\newblock In {\em Triangulated categories}, volume 375 of {\em London Math. Soc. Lecture Note Ser.}, pages 76--160. Cambridge Univ. Press, Cambridge, 2010.

\bibitem[KS11]{Keller-Sch:linear-recurrence}
B.~Keller and S.~Scherotzke.
\newblock {Linear recurrence relations for cluster variables of affine quivers}.
\newblock {\em Adv. Math.}, 228(3):1842--1862, 2011.

\bibitem[MG12]{SMG:arith-2-frieze}
S.~Morier-Genoud.
\newblock Arithmetics of 2-friezes.
\newblock {\em J. Algebraic Combin.}, 36(4):515–539, 2012.

\bibitem[MG15]{Sophie-M:survey}
S.~Morier-Genoud.
\newblock {Coxeter's} frieze patterns at the crossroads of algebra, geometry and combinatorics.
\newblock {\em Bull. Lond. Math. Soc.}, 47(6):895--938, 2015.

\bibitem[Mul23]{muller:finite}
G.~Muller.
\newblock {A short proof of the finiteness of Dynkin friezes}.
\newblock arXiv:2312.11394, 2023.

\bibitem[OV90]{vinberg:liegroup}
A.~L. Onishchik and È.~B. Vinberg.
\newblock {\em {Lie groups and algebraic groups}}.
\newblock Springer-Verlag, Berlin, 1990.

\bibitem[Rin12]{Ringel:cluster-additive}
C.~M. Ringel.
\newblock {Cluster-additive functions on stable translation quivers}.
\newblock {\em J. Algebraic Combin.}, 36(3):475--500, 2012.

\bibitem[Rup13]{rupel:greedy-gen-clust-alg}
D.~Rupel.
\newblock Greedy bases in rank 2 generalized cluster algebras.
\newblock arXiv: 1309.2567, 2013.

\bibitem[YZ08]{YZ:Lcc}
S.-W. Yang and A.~Zelevinsky.
\newblock {Cluster algebras of finite type via Coxeter elements and principal minors}.
\newblock {\em Transform. Groups}, 13(3-4):855--895, 2008.

\end{thebibliography}
\end{document}